\documentclass{report}
\usepackage{fullpage,makeidx,epsf,amssymb}
\def\C{{\rm C \kern-.48em\vrule width.06em height.57em depth-.02em \kern.48em}}
\def\sopenC{{\rm C\kern-.15cm\vrule width.6pt height 4.1pt depth-.3pt
 \kern.15cm}}
\def\Z{{{\rm Z}\kern-.28em{\rm Z}}}
\def\N{{{\rm I}\kern-.16em{\rm N}}}
\def\R{{{\rm I}\kern-.16em{\rm R}}}
\def\Re{\mathop{\rm Re}\nolimits}
\def\set#1{\langle #1 \rangle}
\def\sm{\mathop{\cal S}\nolimits}
\def\dm{\mathop{\cal D}\nolimits}
\def\curve{\mathop{\cal C}\nolimits}
\def\sgn{\mathop{\rm sign}\nolimits}
\def\ls{\mathop{\cal L}\nolimits}
\def\order{\mathop{\rm order}\nolimits}
\def\diag{\mathop{\rm diag}\nolimits}
\def\spa{\mathop{\rm span}\nolimits}
\def\belowrightarrow#1{{{{}\over\ #1\ }\kern-1.1em\to}}
\def\mystrut{\vphantom{\vrule width .1em height 1.em depth .7em}}
\def\goback#1{\setbox0\hbox{#1}\kern-\wd0 \relax}
\def\eqbd{\mathop{{:}{=}}}
\def\bdeq{\mathop{{=}{:}}}
\def\y{(1-\lambda)}
\def\nt{\noindent}

\makeindex

\begin{document}

\title{\Huge  Theorems and counterexamples on structured matrices} 
\author{\LARGE Olga V. Holtz \\  
\\
\Large  Department of Mathematics \\
\Large  University of Wisconsin \\  
\Large Madison, Wisconsin 53706 U.S.A. }
\date{August 2000}
\maketitle

\large

\tableofcontents

\chapter*{}

\section{Introduction}   

This thesis is devoted to several problems posed for special classes
of matrices (such as GKK\index{GKK matrix}) and solved using structured matrices (such
as Toeplitz\index{Toeplitz matrix})  belonging to that class. 

The topic of Chapter 1 is GKK\index{GKK matrix} $\tau$-matrices.
\index{tau-matrix@$\tau$-matrix} This notion was
introduced in the 1970's as a response to the Taussky unification
problem posed in the late 50's, which is discussed in Section 1.1.
Section 1.2 is devoted to four conjectures proposed in
the 1970's-1990's on the stability\index{stability} of GKK\index{GKK matrix} 
$\tau$-matrices.
They are all disproved in Section 1.3 using GKK\index{GKK matrix} 
$\tau$-matrices
with additional structure (Toeplitz\index{Toeplitz matrix} and Hessenberg). 
Further properties of the counterexample matrices, which are themselves
not important  for disproving the stability\index{stability} conjectures but seem to be 
worth analyzing, are taken up in Section~1.4. Section~1.5 contains a brief
discussion of open problems related to the topics of earlier sections.

Chapter 2 is centered around the following problem: given a collection
of matrices $(A_\alpha)$ in a special class (such as totally nonnegative\index{totally nonnegative matrix})
bounded in some matrix norm\index{matrix norm} and such that the spectrum\index{spectrum} of $A_\alpha$ lies 
outside a disk of fixed radius with center at zero, determine whether
the collection $(A^{-1}_\alpha)$ is bounded in the same matrix norm\index{matrix norm}. 
For any matrix norm\index{matrix norm}, the answer is yes for matrices of bounded order,
as is shown in Section 2.1. The next sections all deal with collections of
matrices of unbounded order and the `simplest'
$\infty$-norm\index{p@$p$-norm!$\infty$-norm}, the choice also motivated by applications. It is shown 
in Section 2.2 that the answer is still yes for totally nonnegative\index{totally nonnegative matrix} 
Hermitian\index{Hermitian matrix} matrices. However, the answer is no for positive definite\index{Hermitian matrix!positive definite}
Hermitian\index{Hermitian matrix} matrices. Section 2.3 contains a pertinent 
counterexample and a variation of it both based on the Hilbert matrix. 
A counterexample for the class of Hermitian\index{Hermitian matrix} Toeplitz\index{Toeplitz matrix} matrices 
is obtained in Section~2.4. Finally, an interesting question of 
the same type arising in spline theory is discussed in Section~2.5.   

\section{Notation}

The following conventions are used throughout the thesis.

To make a clear distinction between equality and equality by definition,
the latter is denoted by $\eqbd$. The symbol $\#$ denotes the cardinality
of a set. The set $\{1, \ldots, n\}$ for $n\in \N$ is denoted by $\set{n}$. 
For $p$, $q\in \Z$, let $$ p{:}q \; \eqbd\cases{ \{p,p+1,\ldots,q-1, q\} & 
if $p\leq q$ \cr \emptyset &
                 otherwise}. $$
For any $x\in \R$, set
                 $$ x_+\eqbd\cases{ x & if $x>0$ \cr 0 &
                 otherwise \cr}.$$

The linear space of (column-)vectors with $n$ entries in $\C$ is
denoted by $\C^n$. The $j$th vector of the standard
basis of $\C^n$, i.e., the vector with $1$ at the $j$th position
and zeros elsewhere\index{unit vector}, will be denoted by $e_j$.
The linear space of all $n\times n$-matrices with entries in $\C$
is denoted by $\C^{n\times n}$.  
If the order of a matrix is not clear from the
context, it will be indicated by the subscript $n\times n$ or simply $n$.
The symbol $I$ stands for the identity matrix of appropriate order.
Both the zero matrix and the number zero are denoted by $0$.  

\section{Basic matrix notions}

Given a matrix $A\in \C^{n\times n}$,  let $A(\alpha,\beta)$ denote the
submatrix of $A$ whose rows are indexed by $\alpha$ and columns by
$\beta$ ($\alpha$, $\beta\in \set{n}$) and let
$A[\alpha,\beta]$ denote $\det A(\alpha,\beta)$ if $\# \alpha=\#
\beta$. For simplicity, $A(\alpha)$ will stand for $A(\alpha,\alpha)$
and $A[\alpha]$ for $A[\alpha,\alpha]$. By definition, 
$A[\emptyset]\eqbd1$. Elements of $A$ are denoted by $a(i,j)$.
A block diagonal matrix with (square) blocks $A$ and $B$ will be denoted by 
$$ \diag(A,B)
(\eqbd \left( \begin{array}{cc} A & 0 \\ 0 & B \end{array} \right)).$$
 
The spectrum\index{spectrum} of $A$, i.e., the multiset of its eigenvalues
(with each eigenvalue repeated according to its multiplicity),
is denoted by $\sigma(A)$. The spectral radius\index{spectral radius} of $A$
is denoted by $\varrho(A)(\eqbd \max|\sigma(A)|)$. 

The inequality $A\geq 0$ ($>0$) means that $A$ is entrywise nonnegative
(positive). $A\geq B$ ($A>B$) means, by definition, that $A-B\geq 0$
($A-B>0$).

 A norm $\| \cdot \|$ on the space $\C^{n\times n}$ is a {\em matrix norm\index{matrix norm}\/} 
if it satisfies the inequality
$$ \| AB \| \leq \|A\| \cdot \|B\| \qquad \forall \; A, B\in \C^{n\times n}.$$
A matrix norm\index{matrix norm} $\| \cdot \|_o$ is the {\em operator norm\index{operator norm}\/} subordinate to
the norm $\|\cdot \|$ on $\C^n$ if
$$ \|A\|= \sup_{v\in \sopenC^n \setminus \{0\}} {\|Av\|\over \|v\|} \qquad \forall \; 
A\in \C^{n\times n}. $$ In particular, the {\em $p$-norm\index{p@$p$-norm}\/} $\|\cdot\|_p$ 
($1\leq p \leq \infty$) on $\C^{n\times n}$ is the operator norm\index{operator norm} 
subordinate to the $p$-norm\index{p@$p$-norm} 
$$ \|v\|_p\eqbd \cases{ \left( \sum_{i=1}^n |v(i)|^p \right)^{1/p} & if 
$1\leq p <\infty$ \cr \max_{i=1,\ldots,n} |v(i)| & if $p=\infty$} $$
on the space $\C^n$. The {\em condition number\/}\index{condition number} 
of an invertible matrix $A$ (for the norm $\| \cdot\|$) is 
the product $\|A\|\cdot \|A^{-1}\|$.

A matrix $T\in \C^{n\times n}$ is {\em Toeplitz\index{Toeplitz matrix}\/} if it has the form 
$T\bdeq (\tau(i-j))_{i,j=1}^n$ for some $(\tau(i))_{i=1-n}^{n-1}$.
A matrix is {\em Hessenberg\/} if the entries on its first subdiagonal
are all equal to $1$ and the entries below that subdiagonal are zero. 

\chapter{Eigenvalues of GKK\index{GKK matrix} matrices} 

\section{Taussky unification problem}

A matrix $A$ is called {\em totally nonnegative\index{totally nonnegative matrix}\/}
if $A[\alpha, \beta]\geq 0$  for all $\alpha$, $\beta\in \set{n}$
with $\# \alpha=\# \beta$. $A$ is called an {\em $M$-matrix\index{m@$M$-matrix}\/} if 
$A=r I-P$ where $P\geq 0$ and $r>\varrho(P)$\index{spectral radius}. 
For more than a dozen
other ways to define $M$-matrices\index{m@$M$-matrix}, see~\cite{BP}.
 
A matrix $A$ is called a {\em $P$-matrix\index{p@$P$-matrix} \/} if $A[\alpha]>0$
for all $\alpha\subseteq \set{n}$. $A$ is said to be {\em sign-symmetric\index{sign-symmetric}\/}
if $$A[\alpha,\beta] A[\beta,\alpha]\geq 0 \qquad \forall \alpha,\beta\in
\set{n}, \quad \# \alpha=\# \beta.$$
$A$ is called {\em weakly sign-symmetric\index{sign-symmetric!weakly}\/} if
$$ A[\alpha,\beta] A[\beta,\alpha]\geq 0$$
for all $ \qquad \forall \alpha,\beta\in
\set{n}$, with
\begin{equation}
 \# \alpha=\# \beta=\# \alpha \cup \beta -1. \label{almost} 
\end{equation}
The minors $A[\alpha,\beta]$ with the property~(\ref{almost}) 
are sometimes called  {\em almost principal\/}.

Weakly sign-symmetric\index{sign-symmetric!weakly} $P$-matrices\index{p@$P$-matrix} are also called {$\!$\em GKK\index{GKK matrix}\/} 
after Gantmacher, Krein, and Kotelyansky. 

Let $$l(A)\eqbd \min \sigma(A) \cap \R,$$ with the understanding that, in this setting, $\min \emptyset=\infty$. A matrix $A$ is
called an {\em $\omega$-matrix\/} if it has eigenvalue\index{eigenvalue monotonicity} monotonicity in the sense that
$$ l(A(\alpha,\alpha))\leq l(A(\beta,\beta))<\infty \qquad \mbox{whenever}
\quad \emptyset \neq \beta \subseteq \alpha \subseteq \set{n}.$$ $A$ is a
{\em $\tau$-\index{tau-matrix@$\tau$-matrix}matrix\/} if, in addition, $l(A) \geq 0$. 

Hermitian\index{Hermitian matrix} positive definite\index{Hermitian matrix!positive definite}, nonsingular totally nonnegative\index{totally nonnegative matrix}, and $M$-matrices\index{m@$M$-matrix} 
all enjoy  positivity of principal minors, weak sign symmetry, and 
eigenvalue\index{eigenvalue monotonicity} monotonicity. In fact, these properties were singled out 
as a response to the `unification problem' for the above-mentioned  
three classes of matrices that stems from a research problem posed 
by O.~Taussky~\cite{TT}.\footnote{O. Taussky pointed out in~\cite{TT} that 
similar  theorems were known for some positive matrices and for positive
definite Hermitian matrices, for which the then available proofs were
different, and asked for a unified treatment of both cases. She gave
four examples of such similar theorems, two of which illustrate common 
properties of totally nonnegative and positive definite Hermitian matrices. 
The term `Taussky unification problem' was later taken to mean a much wider 
class of problems.}
 
The fact that Hermitian\index{Hermitian matrix} positive definite\index{Hermitian matrix!positive definite} matrices are sign-symmetric\index{sign-symmetric}
$P$-matrices\index{p@$P$-matrix} (a property stronger than being GKK\index{GKK matrix} matrices) is 
standard. The eigenvalue interlacing\index{eigenvalue interlacing} property of Hermitian\index{Hermitian matrix} 
matrices~(\cite{Ca} or, e.g.,~\cite[p.59]{Bh}) implies their eigenvalue\index{eigenvalue monotonicity} 
monotonicity.
Directly from the definition, nonsingular totally nonnegative\index{totally nonnegative matrix} 
matrices are sign-symmetric\index{sign-symmetric} with nonnegative principal minors. Their 
eigenvalue\index{eigenvalue monotonicity} monotonicity was proved by Friedland~\cite{Fr}. This fact and 
the spectral theory of totally nonnegative\index{totally nonnegative matrix} matrices show that all principal
minors of a nonsingular totally nonnegative\index{totally nonnegative matrix} matrix are in fact positive.

The Perron-Frobenius\index{Perron-Frobenius theory} spectral theory of
nonnegative matrices shows that $M$-matrices\index{m@$M$-matrix} are $P$- and $\omega$-matrices.
The weak sign symmetry of $M$-matrices\index{m@$M$-matrix} was proved by Carlson~\cite{C1}.

\section{Stability conjectures on GKK $\tau$-matrices} 

Yet another property shared by Hermitian\index{Hermitian matrix} positive definite\index{Hermitian matrix!positive definite}, totally
nonnegative, and $M$-matrices\index{m@$M$-matrix} is their positive stability\index{stability}. To recall,
a matrix is called {\em positive (negative) stable\index{stability}\/} if its spectrum\index{spectrum}
lies entirely in the open right (left) half plane.  In the sequel, the 
term `positive stable\index{stability}' will be usually shortened to simply `{\em stable\index{stability}}'. 

Hermitian\index{Hermitian matrix} positive definite\index{Hermitian matrix!positive definite} and totally nonnegative\index{totally nonnegative matrix} matrices are obviously
stable\index{stability} (having only positive eigenvalues), while the
stability\index{stability} of 
$M$-matrices\index{m@$M$-matrix} follows from the Perron-Frobenius theory.

The natural question arising from this observation was whether some 
combination of the properties from Section~1.1, viz., positivity 
of principal minors, weak sign-symmetry, and  eigenvalue\index{eigenvalue monotonicity} monotonicity,
is sufficient to guarantee stability\index{stability}. (None of those properties alone is
sufficient, which can be checked by simple examples.)

Carlson\index{Carlson's conjecture}~\cite{C2} conjectured that the GKK\index{GKK matrix} matrices are stable\index{stability}  and showed his
conjecture to be true for $n\leq 4$.

Engel\index{Engel-Schneider's conjecture} and Schneider~\cite{ES} asked if nonsingular $\tau$-matrices or,
equivalently, {$\omega$-matrices} all whose principal minors are positive 
(see Remark~3.7 in~\cite{ES}), are positive stable\index{stability}. Varga\index{Varga's conjecture}~\cite{V} 
conjectured even more than stability\index{stability}, viz.
$$ |\arg(\lambda - l(A))|\leq \frac{\pi}{2}-\frac{\pi}{n} \qquad
\forall \lambda \in \sigma(A). $$ This inequality was proven for
$n\leq 3$ by Varga (unpublished) as well as by Hershkowitz and 
Berman~\cite{HB} 
and for $n=4$ by Mehrmann~\cite{M}.  In his survey paper~\cite{He},
Hershkowitz\index{Hershkowitz' conjecture} posed the weaker conjecture that $\tau$-\index{tau-matrix@$\tau$-matrix}matrices that are
also GKK\index{GKK matrix} are stable\index{stability}. 

The above conjectures were plausible not only because they were 
verified for matrices of small order, but also due to the following 
two theorems.

The first indicates that there is a certain `forbidden wedge'
around the negative real axis  where eigenvalues of a $P$-matrix\index{p@$P$-matrix} 
cannot lie. (The angle of the wedge depends on the order of the
matrix.) 

\nt {\bf Theorem (Kellogg\index{Kellogg's theorem} \cite{Ke}).\/} {\em 
Let $A\in \C^{n\times n}$ be a $P$-matrix\index{p@$P$-matrix}. Then
$$ |\arg(\lambda)|<\pi-{\pi\over n} \qquad \forall \lambda \in \sigma(A).$$ }

The second  shows that sign symmetry together with positivity of principal
minors is sufficient for stability\index{stability}. 

\nt {\bf Theorem (Carlson\index{Carlson's theorem} \cite{C2}).\/}   {\em 
Sign-symmetric\index{sign-symmetric} $P$-matrices\index{p@$P$-matrix} 
are stable\index{stability}. }

Carlson's elegant proof employs the Cauchy-Binet
\index{Cauchy-Binet formula} formula to show that
 $A^2$ is a $P$-matrix\index{p@$P$-matrix} whenever $A$ is sign-symmetric\index{sign-symmetric}, diagonal stable scaling property\index{diagonal scaling}
 of $P$-matrices\index{p@$P$-matrix} to conclude that $DA$ is 
stable\index{stability} for some
diagonal\index{diagonal matrix} matrix $D$ with positive diagonal, so that 
the homotopy $S(t)\eqbd((1-t)D +tI)A$ preserves sign symmetry as well as
positivity of principal minors, hence the eigenvalues of $S(t)$ cannot
cross the imaginary axis as $t$ runs from $0$ to $1$.

\section{Counterexample}

As is shown in~\cite{Ho}, none of the proposed conjectures is true. 
Described below is a class of GKK\index{GKK matrix} $\tau$-\index{tau-matrix@$\tau$-matrix}matrices which are not even
nonnegative stable, i.e., they do have eigenvalues with negative real part.
This class consists of Toeplitz\index{Toeplitz matrix} 
Hessenberg\index{Hessenberg matrix} matrices $A_{n,k,t}$ of order $n$ that depend on
two more parameters $k\in \N$ and $t\in \R$. In what follows, it will
be shown  that $A_{n,k,t}$  is a GKK\index{GKK matrix} 
$\tau$-\index{tau-matrix@$\tau$-matrix}matrix for any $t\in (0,1)$ 
and that $A_{2k+2,k,t}$ is unstable for sufficiently large $k$ and 
sufficiently small positive $t$.

Let $A$ be an infinite Toeplitz Hessenberg matrix with 
first row $(a_0,a_1,\ldots)$ and let $d_n$ denote its leading 
principal minor of order $n$.
By the Laplace expansion \index{Laplace expansion} by minors of the first row,
\begin{equation}
d_n=\sum_{j=1}^n (-1)^{n-1} a_{j-1} d_{n-j}, \qquad  n\in \N.
\label{basic} 
\end{equation}
This is, in effect, an invertible lower triangular system for the $a_j$.
So, for an arbitrary sequence $(d_1,d_2,\ldots)$, there exists exactly
one Toeplitz Hessenberg matrix having these as its leading principal minors.

With this, let $A_{\infty,k,t}$ be the Toeplitz Hessenberg matrix with leading
principal minors 
$$ d_n=t^{(n-k-1)_+}, \qquad  n\in \N,$$
for some $t\in (0,1)$ and $k\in \N$.
Then equation~(\ref{basic}) becomes  
\begin{equation}
\sum_{j=0}^{n-1} (-1)^ja_j t^{(n-j-k-2)_+} = t^{(n-k-1)_+}, \qquad
n\in \N. \label{a-s}
\end{equation}
Let $A_{n,k,t}$ be the leading principal submatrix of $A_{\infty,k,t}$ 
of order $n$.

Two observations are immediate:
The first $k+1$ entries in the sequence $(d_1,d_2,\ldots)$ equal $1$ and
$d_{k+2}=t$, hence $a_0=1$, $a_j=0$ for $1\leq j \leq k$,  and
$a_{k+1}=(-1)^k(1-t)$. Secondly, the matrix $A_{\infty,k,t}$ is 
Toeplitz\index{Toeplitz matrix}, so all its principal minors indexed by 
$j$ consecutive integers equal $t^{(j-k-1)_+}$.

\subsection{On two characteristics of $A_{\infty,k,t}$}

As is well known, associated to any infinite Toeplitz\index{Toeplitz matrix} matrix
 $T\bdeq(\tau(i-j))_{i,j=0}^\infty$ is its {\em symbol\index{Toeplitz matrix!symbol of}\/} $\sm(T)$, i.e., 
the Laurent\index{Laurent series} series
$$ \sm(T)(s)\eqbd \sum_{j=-\infty}^\infty \tau(-j)s^j. $$  It is also 
convenient to introduce the Taylor\index{Taylor series}  series $\dm(T)$
$$ \dm(T)(s,\lambda)\eqbd \sum_{j=0}^\infty D_j(\lambda)s^j$$
where 
$$ D_j(\lambda) \eqbd \cases{ 1 & if $j=0$ \cr 
\det (T(1{:}j)-\lambda I_j) & if $j\in \N$}.$$
To avoid cumbersome notation, $D_j(\lambda)$ will be used in the sequel,
but one should keep in mind that $D_j(\lambda)$ also depends $T$, i.e., for
$T=A_{\infty,k,t}$, on the 
parameters $k$ and $t$, so should have been denoted by $D_j(\lambda, k,t)$.

\nt {\bf Proposition 1.\/} {\em
\begin{eqnarray}
\sm(A_{\infty,k,t})(s)&=&{1+ts \over s(1+(1-t)\sum_{j=1}^{k+1}(-s)^j)}, \label{sm} \\
\dm(A_{\infty,k,t})(s,\lambda)&=&{1+(1-t)\sum_{j=1}^{k+1}s^j \over 
1+s(\lambda-t)+\lambda(1-t)\sum_{j=2}^{k+2}s^j }. \label{dm}  
\end{eqnarray} }
 
{\em Proof.\/} \hspace{0.3cm} From~(\ref{a-s}), 
\begin{equation}
 a_l=(-1)^l t^{(l-k)_+} +\sum_{j=1}^{l}(-1)^{l-j}t^{(l-j-k)_+}a_{j-1} \label{replace}
\end{equation}
for all $l\in \N$. Let $\Phi(s)\eqbd\sm(A_{\infty,k,t})(s)-1/s$.
Replacing each $a_j$, $j\in \N$, by its expansion~(\ref{replace}) 
and collecting terms with the factor $t^{(j-k)_+}$ for each $j$, one obtains
$$ \Phi(s)= a_0 -\Phi(s)\sum_{j=1}^\infty t^{(j-k-1)_+}
(-s)^{j} +\sum_{j=1}^\infty t^{(j-k)_+} (-s)^j.$$
Recall that $a_0=1$ and compute the sums in the right-hand side.
This yields
\begin{equation}
 {1+ts-(1-t)(-s)^{k+2} \over (1+s)(1+ts) } \Phi(s) = {1+ts-(1-t)(s)^{k+1}
 \over  (1+s)(1+ts)}, \label{Phi} \end{equation}
hence $$
 \sm(A_{n,k,t})(s)=\Phi(s) +{1\over s}= {(1+s)(1+ts) \over
 s(1+ts-(1-t)(-s)^{k+2})}, $$ 
 which, after canceling $(1+s)$, gives~(\ref{sm}).

Now observe that, similarly to~(\ref{a-s}), expansion \index{Laplace expansion}
 of the determinant
$D_n(\lambda)$ by its first row gives 
\begin{equation}
(a_0-\lambda) D_{n-1}+\sum_{j=1}^{n-1} (-1)^ja_j D_j(\lambda) = D_n(\lambda), 
\qquad n\in \N. \label{d-s}
\end{equation}
That is, the coefficient of $s^j$ in $\dm(A_{\infty,k,t})(s,\lambda)$ equals 
the coefficient of $s^{j-1}$ in the product of $\dm(A_{\infty,k,t})(s,
\lambda)$ and $\Phi(-s)-\lambda$
whenever $j\geq 1$. With the zeroth coefficient taken into account, this means
$$ \dm(A_{\infty,k,t})(s,\lambda)= s(\Phi(-s)-\lambda)\dm(A_{\infty,k,t})(s,\lambda)
+1,$$ which, together with~(\ref{Phi}), implies~(\ref{dm}). \hfill $\Box$

\nt {\bf  Corollary.\/} {\em The polynomials $D_j$ satisfy the recurrence relation
\begin{equation}
D_j(\lambda)+(\lambda-t)D_{j-1}(\lambda)+\lambda(1-t)\sum_{l=2}^{k+2}
D_{j-l}(\lambda)=0 \qquad \forall j\geq k+2. 
\label{d-recur} 
\end{equation} }

{\em Proof.\/} \hfill Compare the coefficient of $s^j$ for $j\geq k+2$ in the
numerator (zero) with that in the product of $\dm(A_{\infty,k,t})$ and 
the denominator in~(\ref{dm}). \hfill   $\Box$

These results will be revisited in subsection 1.3.3.  Let us now show that
$A_{\infty,k,t}$ is GKK\index{GKK matrix}. 

\subsection{$A_{n,k,t}$ are GKK\index{GKK matrix}}

Since $A_{n,k,t}$ is Hessenberg\index{Hessenberg matrix}, the submatrix
$A_{n,k,t}({\set{n} {\setminus} i{:} i{+}j{-}1})$ is block upper
triangular if $1<i\leq i+j-1<n$, so
\begin{equation}
A_{n,k,t}[{\alpha\cup \beta}]=A_{n,k,t}[ \alpha] A_{n,k,t}
[ \beta] \quad {\rm whenever} \quad i<j-1 \;\; \hbox{\rm for all
} i\in \alpha, \; j \in \beta. \label{*}
\end{equation}
This shows that $A_{n,k,t}$ is a $P$-matrix\index{p@$P$-matrix}.

Fortunately, there is no need to verify the weak sign symmetry of
a $P$-matrix\index{p@$P$-matrix} by computing its almost principal minors. Instead, one
makes use of the following remarkable fact.

\nt {\bf Theorem (Gantmacher, Krein~\cite[p.55]{GK}, and Carlson\index{Gantmacher-Krein-Carlson's theorem}~\cite{C1}).\/} {\em
A $P$-matrix\index{p@$P$-matrix} $A$ is GKK\index{GKK matrix} iff it satisfies the generalized Hadamard-Fisher\index{Hadamard-Fisher inequality} 
inequality
\begin{equation}
A[\alpha] A[\beta] \geq A[{\alpha\cup \beta}]
A[{\alpha \cap \beta}] \qquad \forall \alpha,\beta \subseteq
\set{n}. \label{HF}
\end{equation} }

Since $0<t<1$ and 
$$ (x+y)_+ + (x+z)_+ \leq x_+ + (x+y+z)_+ \qquad
\forall x \quad \forall y,z\geq 0, $$ one obtains
\begin{equation}
\begin{array}{l}
\mystrut A_{n,k,t}[{i{:} i{+}j{-}1}] \cdot A_{n,k,t}[{l{:}
l{+}m{-}1}]\\ = \quad \mystrut t^{(j-k-1)_+ + (m-k-1)_+}\geq
t^{(l+m-i-k-1)_+ + (i+j-l-k-1)_+}\\ =\quad \mystrut 
A_{n,k,t}[{i{:} l{+}m{-}1}]\cdot A_{n,k,t}[{l{:} i{+}j{-}1} \\
\end{array} \qquad {\rm if} \; l\leq i{+}j{-}1. \label{**} 
\end{equation}
Together with~(\ref{*}), ~(\ref{**}) shows that $A_{n,k,t}$
satisfies~(\ref{HF}) if $\alpha$, $\beta$ are sets of consecutive
integers.

To prove~(\ref{HF}) in general, first make a definition. Call the
subsets $\alpha$, $\beta \subseteq \set{n}$ {\em separated\/} if
$|p-q|>1$ $\forall p\in \alpha$, $q\in \beta$. Suppose $\alpha$,
$\beta_1$, $\ldots$, $\beta_j \subseteq \set{n}$ are sets of
consecutive integers, $\beta_i$ ($i=1,\ldots, j$) are pairwise
separated, and
\begin{equation} 
 \mbox{ for any } i=1,\ldots, j, \mbox{ there exist } p\in \beta_i
 \mbox{ and } q\in \alpha  \mbox{ such that } |p-q|\leq 1. 
\label{suppose}
\end{equation}  
Then $A_{n,k,t}$, $\alpha$, and $\beta\eqbd\cup_{i=1}^j \beta_i$
satisfy~(\ref{HF}). Indeed,~(\ref{HF}) holds for $\alpha$ and
$\beta_1$.  If $1\leq l<j$, then, assuming~(\ref{HF}) for $\alpha$ and
$\gamma_l\eqbd \cup_{i=1}^l \beta_i$, 
\begin{eqnarray*}
&&  A_{n,k,t} [{\alpha}] A_{n,k,t} [{\gamma_{l+1}}]
=A_{n,k,t} [{\alpha}] A_{n,k,t} [{\gamma_l}] A_{n,k,t}
[{\beta_{l+1}}] \\ && \quad \geq  A_{n,k,t} [{\alpha\cup \gamma_l}]
A_{n,k,t} [ {\alpha \cap \gamma_l}] A_{n,k,t}[{\beta_{l+1}}].
\end{eqnarray*}
Due to~(\ref{suppose}), $\alpha \cup \gamma_l$ is a set of consecutive 
integers, so an application of~(\ref{HF}) yields
$$ A_{n,k,t} [{\alpha \cup \gamma_l}] A_{n,k,t}[{\beta_{l+1}}]
\geq A_{n,k,t} [{\alpha\cup \gamma_{l+1}}] A_{n,k,t} [{
(\alpha \cup \gamma_{l})\cap\beta_{l+1}}]. $$
But $(\alpha \cup \gamma_l)\cap \beta_{l+1}=\alpha\cap \beta_{l+1}$
since the sets $\beta_i$ are pairwise disjoint. So, 
\begin{eqnarray} 
 A_{n,k,t}[{\alpha}] A_{n,k,t}[{\gamma_{l+1}}] 
& \geq &  A_{n,k,t}[{\alpha \cup \gamma_{l+1}}] A_{n,k,t}[{\alpha
\cap \gamma_l}] A_{n,k,t}[{\alpha \cap \beta_{l+1}}] \nonumber \\
&=& A_{n,k,t} [{\alpha \cup \gamma_{l+1}}] A_{n,k,t}  
[{\alpha\cap \gamma_{l+1}}]. \nonumber
\end{eqnarray}
Now, given a set of consecutive integers $\alpha\subseteq \set{n}$ and
any set $\beta \subseteq \set{n}$, write $\beta=\gamma_1\cup
\gamma_2$ where $\gamma_1\eqbd\cup_{i=1}^l\beta_i$,
$\gamma_2\eqbd\cup_{i=l+1}^{l+m} \beta_i$, all $\beta_i$ ($i=1,\ldots,
l+m$) are separated, and $\beta_i$ satisfies~(\ref{suppose}) if and only if $i\leq l$. Then
\begin{eqnarray*}
&& A_{n,k,t} [{\alpha}] A_{n,k,t}[{\beta}] =
A_{n,k,t} [{\alpha}] A_{n,k,t}[{\gamma_1}] A_{n,k,t}
[{\gamma_2}] \geq  A_{n,k,t} [{\alpha \cup \gamma_1}] A_{n,k,t}[{\alpha
\cap \gamma_1}] A_{n,k,t}[{\gamma_2}] \\
&& \quad = A_{n,k,t} [{\alpha \cup \gamma_1 \cup \gamma_2}] 
A_{n,k,t}  [{\alpha \cap \gamma_1}]  = A_{n,k,t} [{\alpha \cup \beta}] 
A_{n,k,t} [{\alpha \cap \beta}].
\end{eqnarray*}   
In other words, $A_{n,k,t}$ satisfies~(\ref{HF}) if $\alpha\subseteq \set{n}$
is a set of consecutive integers and $\beta \subseteq \set{n}$ is arbitrary.

Finally, if $\alpha_1$, $\alpha_2$, $\beta\subseteq \set{n}$, the sets
$\alpha_i$ ($i=1,2$) are separated, ~(\ref{HF}) holds for $\alpha_1$
and $\beta$, and $\alpha_2$ is a set of consecutive integers, 
then~(\ref{HF}) holds for $\alpha\eqbd\alpha_1\cup \alpha_2$ and $\beta$:
\begin{eqnarray*}
&& A_{n,k,t}[{\alpha}] A_{n,k,t}[{\beta}] =
A_{n,k,t}[{\alpha_1}] A_{n,k,t}[{\alpha_2}]
A_{n,k,t}[{\beta}] \\ & & \quad \geq  A_{n,k,t} [{\alpha_1 \cup
\beta}] A_{n,k,t} [{\alpha_1\cap \beta}] A_{n,k,t}
[{\alpha_2}]  \\ && \quad \geq A_{n,k,t} [{(\alpha_1 \cup
\beta) \cup \alpha_2}] A_{n,k,t} [{ (\alpha_1\cup \beta) \cap
\alpha_2}]  A_{n,k,t} [{\alpha_1 \cap \beta}] \\
 && \quad = A_{n,k,t} [{\alpha \cup \beta}] A_{n,k,t} [{\alpha_1\cap
\beta}] A_{n,k,t} [{\alpha_2 \cap \beta}] \\ && \quad =A_{n,k,t}
[{\alpha \cup \beta}] A_{n,k,t} [{\alpha \cap \beta}].
\end{eqnarray*}
So, by induction on the number of `components' of $\alpha$,~(\ref{HF}) 
holds for any $\alpha$, $\beta\subseteq \set{n}$. Thus, an application of
the Gantmacher-Krein-Carlson\index{Gantmacher-Krein-Carlson's theorem} Theorem concludes the proof of the following.

\nt {\bf Proposition 2.\/} {\em The matrices $A_{n,k,t}$ are GKK\index{GKK matrix} 
for any $n$, $k\in \N$ and $t\in (0,1)$.} \hfill $\Box$ 

\subsection{$A_{n,k,t}$ are $\tau$-matrices}

\index{tau-matrix@$\tau$-matrix}
Now check that $A_{n,k,t}$ have eigenvalue\index{eigenvalue monotonicity} monotonicity for any $n\in \N$
and $t\in(0,1)$.
If $\set{n}\supseteq \alpha=\cup_{i=1}^j \alpha_i$ is the union of 
separated sets of consecutive integers, then 
$$ \det (A_{n,k,t}({\alpha})
-\lambda I)= \prod_{i=1}^j \det (A({\alpha_i})-\lambda I) $$ 
since $A_{n,k,t}-\lambda I$ is Hessenberg\index{Hessenberg matrix} (the same observation earlier 
led to~(\ref{*})). Since $A_{n,k,t}-\lambda I$ is Toeplitz\index{Toeplitz matrix}, the  product
in the right hand side equals $\prod_{i=1}^j D_{\set{\# \alpha_i}}(\lambda)$. 
Hence, to prove eigenvalue\index{eigenvalue monotonicity} monotonicity 
of $A_{n,k,t}$ it is enough to prove it for
leading principal submatrices of $A_{n,k,t}$ only, i.e., to show
$$ l(A_{n+1,k,t})\leq l(A_{n,k,t}) \qquad \forall n\in \N.$$

First recall that a $P$-matrix\index{p@$P$-matrix} has no nonpositive real eigenvalues, since
the coefficients of its characteristic polynomial strictly alternate in sign
(Kellogg\index{Kellogg's theorem}'s theorem gives a stronger statement, but one does not need it here).
So, $l(A_{n,k,t})>0$ for all $n\in \N$.
Now note that $D_n(\lambda)=\y^n$ for $n=0,\ldots, k+1$. 
For $n=k+2$, recall that 
$a_{k+1}=(-1)^k (1-t)$, so~(\ref{d-s}) yields $D_n(\lambda)=(1-\lambda)^n-
(1-t)$. Hence,
$$l(A_{k+2,k,t})=1-(1-t)^{1/(k+2)}<1=l(A_{k+1,k,t})=\cdots=l(A_{1,k,t})$$
and $D_j(l(A_{k+2,k,t}))<0$ for all $j<k+2$. 

On the other hand, if
\begin{equation} D_j(l(A_{n,k,t}))<0 \qquad \forall j<n,  \label{sign} 
\end{equation}
then~(\ref{d-recur}) and the positivity of $l(A_{n,k,t})$ imply 
that $$D_{n+1}(l(A_{n,k,t}))=-(1-t)
l(A_{n,k,t})\sum_{j=1}^{k+1}D_{n-j}(l(A_{n,k,t}))<0.$$
But $D_{n+1}(0)=t^{(n-k)_+}>0$, hence $D_{n+1}$ changes its sign
on the interval $(0,l(A_{n,k,t}))$, hence $l(A_{n+1,k,t})<l(A_{n,k,t})$
and~(\ref{sign}) holds for $n+1$. Thus, by induction, the matrices 
$A_{n,k,t}$ have eigenvalue\index{eigenvalue monotonicity} monotonicity. 

\nt {\bf Proposition 3.\/} {\em  $A_{n,k,t}$ is a $\tau$-matrix for any $n,k\in \N$ 
and any $t\in(0,1)$.} \hfill $\Box$

\subsection{$A_{2k+2,k,t}$ is unstable for sufficiently large $k$
and small $t$}

Now let $B_k\eqbd\lim_{t\to 0+} A_{2k+2,k,t}$. The matrix $B_k$ is
Toeplitz\index{Toeplitz matrix} with first column $$(1,1,\underbrace{0, \ldots, 0}_{2k \; {\rm
\scriptstyle terms}})^T$$ and first row $$(1,\underbrace{ 0, \ldots,
0}_{k\; {\rm \scriptstyle terms}}, (-1)^k, (-1)^k,
\underbrace{0,\ldots, 0}_{k-1 \; {\rm \scriptstyle terms}}).$$
Show that there exists $K\in \N$ such that, for all $k>K$, $B_k$ has
an eigenvalue $\lambda$ with $\Re \lambda<0$.  As the eigenvalues
depend continuously on the entries of the matrix, this will
demonstrate that, for any $k>K$, there exists $t\in (0,1)$
such that the GKK\index{GKK matrix} $\tau$-\index{tau-matrix@$\tau$-matrix}matrix $A_{2k+2,k,t}$ has an
eigenvalue with negative real part.

The polynomial $D_{2k+2}$ has a root with negative real part iff
the polynomial $\psi_k$ defined by
$$\psi_k(\lambda)\eqbd\frac{D_{2k+2}(-\lambda)}{(1+\lambda)^{k-1}}=
(1+\lambda)^{k+3}-(k+1)(1+\lambda)+k$$ has a root with positive real
part. Since
$$\psi_k(\lambda)=\lambda \biggl[ \sum_{j=0}^{k+1} {k+3 \choose j}
\lambda^{k+3-j-1}+2 \biggr],$$ it is, in turn, enough to show that
$\eta_k$ defined by
$$\eta_k(\lambda)\eqbd\lambda^{k+3}\psi_k\left(\frac{1}{\lambda}\right)=
2\lambda^{k+2}+\sum_{j=2}^{k+3} {k+3 \choose j} \lambda^{k+3-j}$$ has
a root with positive real part. 

One is now in the position to apply the classical negative stability\index{stability}
criterion of Hurwitz to the polynomial $\eta_k$. It is
more efficient, however, to use the following condition necessary
for nonpositive stability\index{stability} due to Ando\index{Ando's theorem}.

\nt{\bf Theorem (Ando\index{Ando's theorem} \cite{An}).\/} 
{\em The Hurwitz\index{Hurwitz matrix} matrix 
$$ \left( \begin{array}{ccccc}
    a_1 & a_3 & a_5 & a_7 & \cdots \\
    a_0 & a_2 & a_4 & a_6 & \cdots \\
    0 & a_1 & a_3 & a_5   & \cdots \\
    0 & a_0 & a_2 & a_4   & \cdots \\ 
    \vdots & \vdots & \vdots & \vdots & \ddots 
   \end{array} \right)_{d\times d} $$
of a nonpositive stable\index{stability} polynomial $f(x)\bdeq
\sum_{j=0}^d a_j x^{d-j}$ 
of degree $d$ is totally nonnegative\index{totally nonnegative matrix}. }

The Hurwitz\index{Hurwitz matrix} matrix for the polynomial
$\eta_k$ is
$$ H_k\eqbd \left( \begin{array}{cccccc} {k+3 \choose 2} & {k+3
\choose 4} & {k+3 \choose 6} & {k+3 \choose 8} & {k+3 \choose 10} &
\cdots \mystrut \\ 2 & {k+3 \choose 3} & {k+3 \choose 5} & {k+3
\choose 7} & {k+3 \choose 9} & \cdots \mystrut \\ 0 & {k+3 \choose 2}
& {k+3 \choose 4} & {k+3 \choose 6} & {k+3 \choose 8} & \cdots
\mystrut \\ 0 & 2 & {k+3 \choose 3} & {k+3 \choose 5} & {k+3 \choose
7} & \cdots \mystrut \\ 0 & 0 & {k+3 \choose 2} & {k+3 \choose 4} &
{k+3 \choose 6} & \cdots \mystrut \\ \vdots & \vdots & \vdots & \vdots
& \vdots & \ddots \mystrut
\end{array} \right)_{(k+2)\times (k+2)} \goback{$\scriptstyle (k+2)(k+2)$}. $$
Compute the minor $H_k[{2{:}5}]$, taking out the factors $k+3 \choose 2$, 
$k+3 \choose 4$, $k+3 \choose 6$ from its second, third, and fourth columns
respectively. This gives
$$ H_k[{2{:}5}]=-\frac{1}{132300}(3k^3-49k^2-210k-318)(k+4)^2(k+5)
{k+3 \choose 2}{k+3 \choose 4} {k+3\choose 6}.$$ Thus, 
$H_k[{2{:}5}]<0$ for large enough $k$, precisely, for all
$k>20$. So, for $k>20$, $\eta_k$ has a
zero with positive real part, therefore, $D_{2k+2}$ has a zero 
with negative real part. 

This completes the proof of the following.\index{tau-matrix@$\tau$-matrix}

\nt {\bf Theorem 1.\/} {\em The GKK\index{GKK matrix} $\tau$-matrices
$A_{2k+2,k,t}$ are unstable for sufficiently large $k$ and sufficiently
small positive $t$.} \hfill $\Box$ 

\subsection{Numerics} 
To illustrate the result, consider the matrix $A_{44,21,1/2}$, i.e., 
the Toeplitz\index{Toeplitz matrix} matrix whose first column is $$(1,1, \underbrace{0,\ldots,
0}_{42 \; {\rm \scriptstyle terms}})^T$$ 
and first row is
$$ (1, \underbrace{0,\ldots, 0}_{21 \; {\rm \scriptstyle terms}}, -1/2,
-1/2^2, 1/2^3, -1/2^4, \ldots, -1/2^{22})$$
 and the limit matrix
$B_{21}$, with the same first column as $A_{44,21,1/2}$ and first 
row equal to $$(1, \underbrace{0,\ldots,0}_{21 \; 
{\rm \scriptstyle terms}}, -1, -1, \underbrace{0, \ldots, 
0}_{20 \; {\rm \scriptstyle terms}}).$$
According to MATLAB,
the two eigenvalues with minimal real part of the first matrix are
$$-2.809929189497896\cdot 10^{-2} \pm 3.275076252367531\cdot 10^{-1}i;$$ 
those of the second are $$-3.420708309454068\cdot 10^{-2} \pm
3.400425852703498\cdot 10^{-1}i.$$ Further MATLAB calculations suggest
that the 
matrices $A_{28,13,t}$ are already unstable for sufficiently small $t$.

\section{More on the counterexample matrices} 

\subsection{The sign pattern of $A^{-1}_{n,k,t}$}
The matrices $A_{n,k,t}$ turn out to be quite `close' to being
sign-symmetric\index{sign-symmetric} 
matrices for ${n\leq 2k+2}$. Namely, all their minors of order $n-1$ are nonnegative. Since $\det 
A_{n,k,t}>0$ whenever $t> 0$, this property can be restated as follows.

\nt {\bf Proposition 4.\/} {\em $A_{n,k,t}^{-1}$ is 
checkerboard\index{checkerboard matrix} for $n\leq 2k+2$. }

The proof will make use of the Gohberg-Semencul formula for inverting
a Toeplitz\index{Toeplitz matrix} matrix.

\nt {\bf Theorem (Gohberg and Semencul \cite{GS} or, e.g., \cite[p.21]{HR}).\/}
{\em If the Toeplitz\index{Toeplitz matrix} matrix $T$ of order $n$ is 
invertible, then  
\begin{eqnarray}
T^{-1}& = &x_0^{-1} \bigg\{ \left( \begin{array}{cccc}
x_0 & 0 & \cdots & 0 \\ x_1 & x_0 & \cdots & 0 \\ \vdots &  
\vdots & \ddots & \vdots \\ x_{n-1} & x_{n-2} & \cdots & x_0 \end{array}
\right) \left( \begin{array}{cccc} y_0 & y_{-1} & \cdots & y_{-n+1} \\
0 & y_0 & \cdots & y_{-n+2} \\ \vdots & \vdots & \ddots & \vdots \\
0 & 0 & \cdots & y_0 \end{array} \right) - \nonumber \\
&&  \left(\begin{array}{cccccc} 0 & 0 &  \cdots & 0 & 0 \\
y_{-n+1} & 0 &  \cdots & 0 & 0 \\ y_{-n+2} & y_{-n+1} &  
\cdots & 0 & 0 \\ \vdots & \vdots &  \ddots & \vdots & \vdots \\
y_{-1} & y_{-2} &  \cdots & y_{-n+1} & 0 \end{array} \right)
 \left(\begin{array}{ccccc} 0 & x_{n-1} & x_{n-2} & \cdots & x_1 \\
0 & 0 & x_{n-1} & \cdots & x_2 \\ \vdots & \vdots & \vdots & \ddots &
\vdots \\ 0 & 0 & 0 & \cdots & x_{n-1}  \\
0 & 0 & 0 & \cdots & 0 \end{array} \right) \bigg\}. \label{GS}
\end{eqnarray}
where $x\eqbd(x_0,\ldots,x_{n-1})^T$ is the solution to
$Tx=e_1$ and $y\eqbd(y_{-n+1},\ldots, y_0)^T$ is the solution to $Ty=e_n$
(recall that $e_j$ denotes the $j$-th unit \index{unit vector} vector). } 

{\em Proof of Proposition 4.\/} \hspace{0.3cm} There is nothing to prove
if $n\leq k+1$. The proof for $n\geq k+2$ is by induction 
on $n$. To avoid confusion, let us make the notation more precise.
The goal is to inductively solve equations $A_{n,k,t}X_n=e_1^{(n)}$
and $A_{n,k,t}Y_n=e_n^{(n)}$.
 
Let $s\eqbd1/t$. Prove by induction that \begin{eqnarray}
x_l&=&\cases{ (-1)^l s^{l+1} & if $l\leq n-k-2$ \cr
(-1)^l s^{n-k-1} & if $l\geq n-k-1$,} \qquad l=0,\ldots, n-1, \label{X} \\
y_l& = & \cases{ s & if $l=0$ \cr (-1)^{l+1}(1-s) & if $-1\geq l \geq -k-1$ \cr
0 & if $l\leq -k-2$, } \qquad l=-n+1,\ldots,0 \label{Y} 
\end{eqnarray}
for $n\geq k+2$. Indeed, formulas~(\ref{X}) and~(\ref{Y}) hold for 
$n=k+2$ (which can be checked by direct calculation). To justify the
inductive step, one needs to show that
$$ \left( \begin{array}{cc} 1 & B_{n,k,t} \\ e_1^{(n-1)} & A_{n-1,k,t}
\end{array}\right) \left( \begin{array}{c} s \\ -sX_{n-1} \end{array} 
\right) = e_1^{(n)}, \qquad  
\left( \begin{array}{cc} 1 & B_{n,k,t} \\ e_1^{(n-1)} & A_{n-1,k,t} 
\end{array}\right) \left( \begin{array}{c} 0 \\ Y_{n-1} \end{array}
\right) = e_n^{(n)} $$
where $B_{n,k,t}\eqbd A_{n,k,t}(1,2{:} n)$. 
Expanding~(\ref{sm}) to order $2k+2$, one verifies that
$$ a_j = \cases{ (-1)^k(1-t) & if $j=k+1$ \cr (-1)^j
t^{j-k-2}(1-t)^2 & if $k+1<j\leq 2k+2$}. $$
So, \begin{eqnarray*}
B_{n,k,t}X_{n-1}& = & \sum_{i=1}^{n-k-1}a_i(-1)^{k+i-1}s^{n-k-2}=
s^{n-k-2} [ \y-\sum_{i=2}^{n-k-1} t^{i-2}\y^2]=\\
&& s^{n-k-2} [\y-\y^2\frac{1-t^{n-k-2}}{\y}]=s^{n-k-2}\y t^{n-k-2}=\y,
\end{eqnarray*}
so $s(1-B_{n,k,t}X_{n-1})=1$. Since $e_1^{(n-1)}=A_{n-1,k,t}X_{n-1}$
by the inductive hypothesis, this justifies the transition from 
$X_{n-1}$ to $X_n$. Now,
\begin{eqnarray*} 
B_{n,k,t}Y_{n-1}& = & \sum_{i=1}^{n-k-2} a_i(-1)^{n-k-i} (1-s)+a_{n-k-1}s\\
&=& (-1)^{n-1} [(1-s)(\y-\sum_{i=2}^{n-k-2}\y^2 t^{i-2})+t^{n-k-3}\y^2s]\\
&=& (-1)^{n-1} \left[\frac{t-1}{t}(\y-\y^2\frac{1-t^{n-k-3}}{\y})+t^{n-k-3}\y^2
\frac{1}{t} \right]\\ &=& 
(-1)^{n-1}\left[ \frac{t-1}{t}\y t^{n-k-3} +t^{n-k-3}\y^2 \frac{1}{t} \right]
=0 \end{eqnarray*}
and, since $A_{n-1,k,t}Y_{n-1}=e_{n-1}^{(n-1)}$ by the inductive hypothesis,
the inductive step is completed.

This implies \begin{equation} |x_l y_{-m} |\geq |x_{n-m}y_{l-n}|
\qquad \mbox{whenever} \qquad m\geq 1, \; n-m>l.\label{main} \end{equation}
Indeed, if $l-n\leq -k-2$, the right hand side of~(\ref{main}) is
zero. But if $l\geq n-k-1$, then $|x_l|=|x_{n-m}|$ and
$|y_{-m}|=|y_{l-n}|$ (observe that $l-n\neq 0$). 

Since the $(i,j)$ element of the right hand side of~(\ref{GS})
has the form $$\sum_{l=i-\min\{i,j\}}^i x_l y_{i-j-l} -
\sum_{l=n-j}^{n-j+\min\{i,j\}-1} x_l y_{i-j-l} $$
and the terms in both sums have sign $(-1)^{i-j}$, the subtraction
gives $$ \cases{ x_0 y_{i-j}+ (-1)^{i-j} r_{i,j}& 
if $i\leq j$ \cr x_{i-j}y_0 + (-1)^{i-j} r_{i,j}&
if $i\geq j$ } \qquad \mbox{\rm for some} \quad r_{i,j}\geq 0.$$
So, in either case, if the $(i,j)$ element of $A_{n,k,t}^{-1}$ is nonzero,
then its sign is $(-1)^{i-j}$. \hfill $\Box$

\subsection{The spectrum\index{spectrum} of $A_{\infty,k,t}$} 

It is well known (e.g.,~\cite[p.21]{GF}) that the spectrum\index{spectrum!of a Toeplitz matrix} of an infinite 
Toeplitz\index{Toeplitz matrix} 
matrix $T=(\tau(i-j))$ with 
\begin{equation} \sum_{j=-\infty}^\infty |\tau(j)|<\infty \label{summable}
\end{equation}
is the union of the curve
\begin{equation}
\curve(T)=\{ \sm(T)(s) : s\in \C, \;\; |s|=1 \}
\label{curve}
\end{equation}
and those points in $\C\setminus \curve(T)$ whose winding number with respect
to $\curve(T)$ is nonzero. 

\nt {\bf Proposition 5.\/} {\em The spectrum $\sigma(A_{\infty,k,t})$ is the union of
$$ \curve(A_{\infty,k,t})= \{ {1+ts \over s(1+(1-t)\sum_{j=1}^{k+1}(-s)^j)}
: s\in \C, \;\; |s|=1 \}$$
and the set of points in $\C$ whose winding number with respect to 
$\curve(A_{\infty,k,t})$ is nonzero. In particular, $\curve(A_{\infty,k,t})$
contains a negative real point.}

{\em Proof. \/} The value of the symbol $\sm(A_{\infty,k,t})$ (given by~(\ref{sm})) at 
the point $s=-1$ is 
\begin{equation} -(1-t)/(k+2-(k+1)t).  \label{negpoint}
\end{equation} \hfill $\Box$ 
 
Next, let us try to determine the limit set $\ls(A_{\infty,k,t})$ of 
eigenvalues of the finite  sections of $A_{\infty,k,t}$ (the matrices 
$A_{n,k,t}$) as $n$ tends to infinity. It is known (see, e.g.,~\cite{BGS}) 
that $\ls(T)\subseteq \sigma(T)$ for any Toeplitz\index{Toeplitz matrix} matrix 
satisfying~(\ref{summable}). Moreover, if $\sm(T)$ is a rational 
function, there is a characterization of the set $\ls(T)$ 
due to K. M. Day\index{Day's theorem}. 

\nt{\bf  Theorem (Day\index{Day's theorem} \cite{D}).\/}  {\em
Let $T\eqbd(\tau(i-j))_{i,j=0}^\infty$ be a Toeplitz\index{Toeplitz matrix} matrix 
satisfying~(\ref{summable}) with 
symbol\index{Toeplitz matrix!symbol of} $\sm(T)$ that coincides with the expansion of
the rational function  ${F\over GH}$ in the annulus $\{ s\in \C:  R_1<|s|<R_2\}$,
$G$ ($H$) being a polynomial all of whose roots lie in the set 
$|s|\leq R_1$ ($|s|\geq R_2$), $F$ a polynomial having no common factors 
with $GH$, and $ p\eqbd\deg G$. Let  $T_n\eqbd (\tau(i-j))_{i,j=0}^n$. Then
the set 
\begin{equation}
\ls(T) \eqbd\{ \lambda : \lambda=\lim \lambda_m , \;\; \lambda_m\in 
\sigma(T_{i_m})    \} \label{limset}
\end{equation}
coincides with    
$$ \{ \lambda \in \C : |r_p(\lambda)|=|r_{p+1}(\lambda)| \}, $$
where $r_j(\lambda)$ are the roots of the polynomial $R\eqbd F-\lambda GH$
listed in the order of their absolute values $|r_1(\lambda)|\leq 
|r_2(\lambda)| \leq \cdots$.  }

For the problem in hand, $F(s)=1+ts$, $G(s)=s$, $H(s)=1+(1-t)\sum_{j=1}^{k+1}
(-s)^j$, $p=1$, and $R(s)=1+ts-\lambda s (1+(1-t)\sum_{j=1}^{k+1}(-s)^j)$.
It seems hopeless to seek an explicit 
formula for $\ls(A_{\infty,k,t})$ for arbitrary $k$ and $t$.
Even checking whether a particular point belongs to $\ls(A_{\infty,k,t})$
turns out to be rather nontrivial. 

However, one can easily make the following (negative) observation.

\nt {\bf Proposition 6.\/}  {\em The point~(\ref{negpoint}) does not belong to 
$\ls(A_{\infty,k,t})$ for any $k\in \N$ and $t\in (0,1)$. }

{\em Proof.\/} \hspace{0.3cm} Let $\lambda$ equal~(\ref{negpoint}). If $|s| \leq 1$, then 
\begin{equation} |F(s)| \geq 1-t, \label{triangle1} \end{equation} while,
by the triangle inequality,
\begin{equation}
 |\lambda G(s)H(s)| \leq |\lambda| (1+(1-t)(k+1))=1-t, \label{triangle2}
\end{equation}
so $R=0$ iff inequalities~(\ref{triangle1}) and~(\ref{triangle2}) both
become equalities iff $s=-1$. On the other hand, ${d\over ds}R (-1)\neq 0$,
so $s=-1$ is not a double root of $R$. Thus, $|r_1(\lambda)|<|r_2(\lambda)|$
and $\lambda\notin \ls(A_{\infty,k,t})$. \hfill $\Box$

Nevertheless, some eigenvalues of $A_{n,k,t}$ form sequences approaching
points on the negative real axis as $n$ tends to infinity. This can be shown
using the following result of M. Biernacki\index{Biernacki's theorem}.

\nt {\bf Theorem (Biernacki~\cite{Bi} as cited in~\cite{ScS}).\/} {\em 
Let $p$, $q\in \N$ be relatively prime and let 
$f(s)\eqbd s^{-p}+s^q -\lambda$.
Then the set $$ \{ \lambda : |r_p(\lambda)|=|r_{p+1}(\lambda)| \}$$
is the star-shaped curve
$$ S:=\{ \lambda=\varepsilon r \}$$
where $\varepsilon$ is any $(p+q)$th root of unity and
$$ 0\leq r \leq (p+q) p^{-p/(p+q)} q^{-q/(p+q)}. $$ }

Here is the more precise formulation of the above claim.

\nt {\bf Proposition 7.\/} {\em For any $k\in \N$, there exists $t(k)\in (0,1)$
such that the set $\ls(A_{\infty,k,t})$ contains a segment on the negative 
real axis for any $t\in (0,t(k))$.}

{\em Proof.\/} \hspace{0.3cm} If $t=0$,  $\widetilde{R}(s)\eqbd
(1+s)R(s)=1+s-\lambda s(1-(-s)^{k+2})$. For $\lambda\neq 0$, the equation $\widetilde{R}(s)=0$
is equivalent to the equation \begin{equation} 
 (-s)^{-1}+(-s)^{k+2}-\mu=0 \qquad {\rm{where}} \qquad 
\mu \eqbd{1-\lambda \over(-\lambda)^{1/(k+3)}}. \end{equation} 
Taking $\varepsilon\eqbd1$ and applying Biernacki\index{Biernacki's theorem}'s theorem
amounts to solving the equation
\begin{equation}
1-\lambda -(-\lambda)^{1/(k+3)}r=0 \label{bier}
\end{equation}
for some $r\in (0,r_k)$ where $r_k\eqbd(k+3)(k+2)^{-1/(k+3)}$. 
Since $r_k>2$ for all $k\geq 0$, the left hand side of~(\ref{bier})
changes sign between $(\lambda=)-1$ and $\lambda=0$ for all 
$r\in(2,r_k)$, therefore, there exists a segment 
$\Lambda\eqbd(\lambda_{min}, \lambda_{max})\subseteq(-1,0)$ such that
$|r_1(\lambda)|=|r_2(\lambda)|$ for the polynomial $\widetilde{R}$
whenever $\lambda\in \Lambda$.    

Without loss of generality, one can assume that $\widetilde{R}$ has no 
double roots when $\lambda\in \Lambda$, since that assumption can rule out
only a finite number of $\lambda$'s. So, one can assume that such
points are already excluded from $\Lambda$.

If, for some $\lambda \in \Lambda$, $\widetilde{R}$ had two distinct 
roots of the form  $s$ and 
$\epsilon |s|$, with $\epsilon \eqbd\pm 1$, 
that would imply that one of the triangle inequalities
\begin{eqnarray*}
 |(\lambda-1)s+\lambda(-s)^{k+3}|&\leq&  |\lambda -1|\cdot |s|+|\lambda|
 \cdot |s|^{k+3} \\
  |(\lambda-1)s+\lambda(-s)^{k+3}|&\geq&  \left|
 |\lambda -1|\cdot |s|-|\lambda| \cdot |s|^{k+3}\right| 
\end{eqnarray*}
becomes equality, which is possible only if $s^{k+2}$ is real, hence, by the 
condition $\tilde{R}(s)=0$, $s$ itself must be real. But the assumption 
that $\widetilde{R}$ has
two roots $s, -s\in \R$ leads to the contradictions $2=0$ if $k$ is even
and $(\lambda-1)s=0$, hence $\lambda=1$, if $k$ is odd. 

So, $\widetilde{R}$ has no roots of the form $s$, $\pm |s|$ for $\lambda\in 
\Lambda$. 

This implies, first of all, that none of the roots of $\widetilde{R}$ 
with minimum absolute value equals $-1$ when $\lambda \in \Lambda$, so
that $r_1(\lambda)$, $r_2(\lambda)$ are also roots of $R$ with minimal
absolute value. Moreover, the set of roots with the smallest absolute
value consists of (possibly, several) distinct non-real conjugate pairs 
$s_j$, $\overline{s_j}$. Since the roots of an algebraic equation
are continuous functions of the coefficients, one of those pairs 
must stay a pair of complex conjugate roots with smallest absolute
value as $t$ runs from $0$ to $t(k)$ for some sufficiently small value
$t(k)$. \hfill $\Box$

To visualize the sets $\sigma(A_{\infty,k,t})$ and
$\ls(A_{\infty,k,t})$, here are four figures, two for $k\eqbd3$, $t\eqbd0.2$
and two for for $k\eqbd10$, $t\eqbd0.4$, drawn by MATLAB.

\begin{center}
\begin{tabular}{c} 
\epsfysize=9.0cm       
\epsfbox{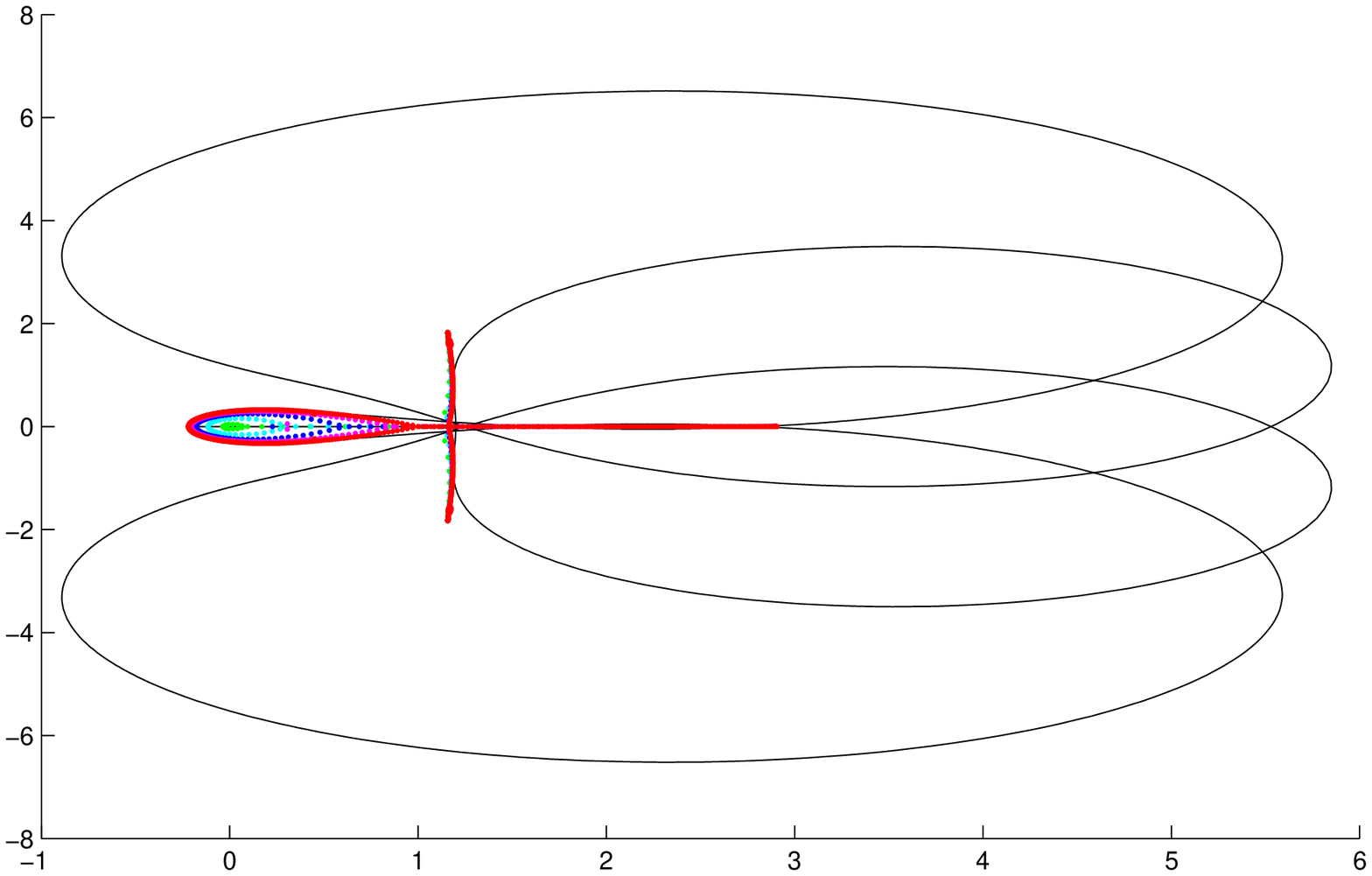} \\
\\
Figure 1. The curve $\curve(A_{\infty,k,t})$ (black) and the sets 
$\sigma(A_{50,k,t})$ (green),  \\
$\sigma(A_{100,k,t})$ (cyan),  $\sigma(A_{200,k,t})$ (blue),
 $\sigma(A_{400,k,t})$ (magenta), \\
$\sigma(A_{800,k,t})$ (red). Here $k=3$, $t=0.2$. \\ 
\end{tabular}

\begin{tabular}{c} 
\epsfysize=9.0cm       
\epsfbox{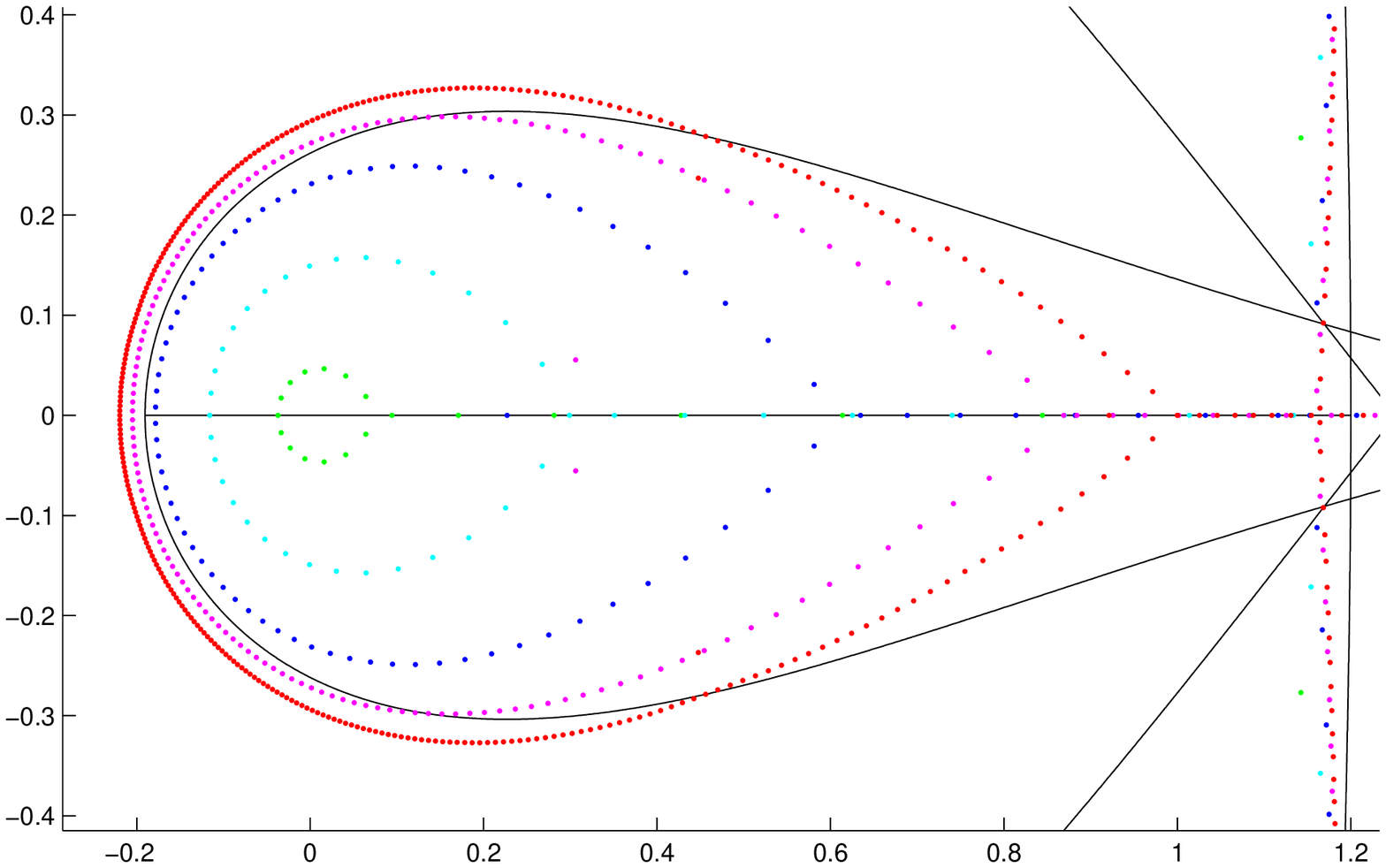} \\
\\
Figure 2. Figure 1 zoomed in the part of $\sigma(A_{\infty,k,t})$ \\
around the negative real axis.\\
\end{tabular}
\end{center}

\begin{center}
\begin{tabular}{c} 
\epsfysize=9.0cm       
\epsfbox{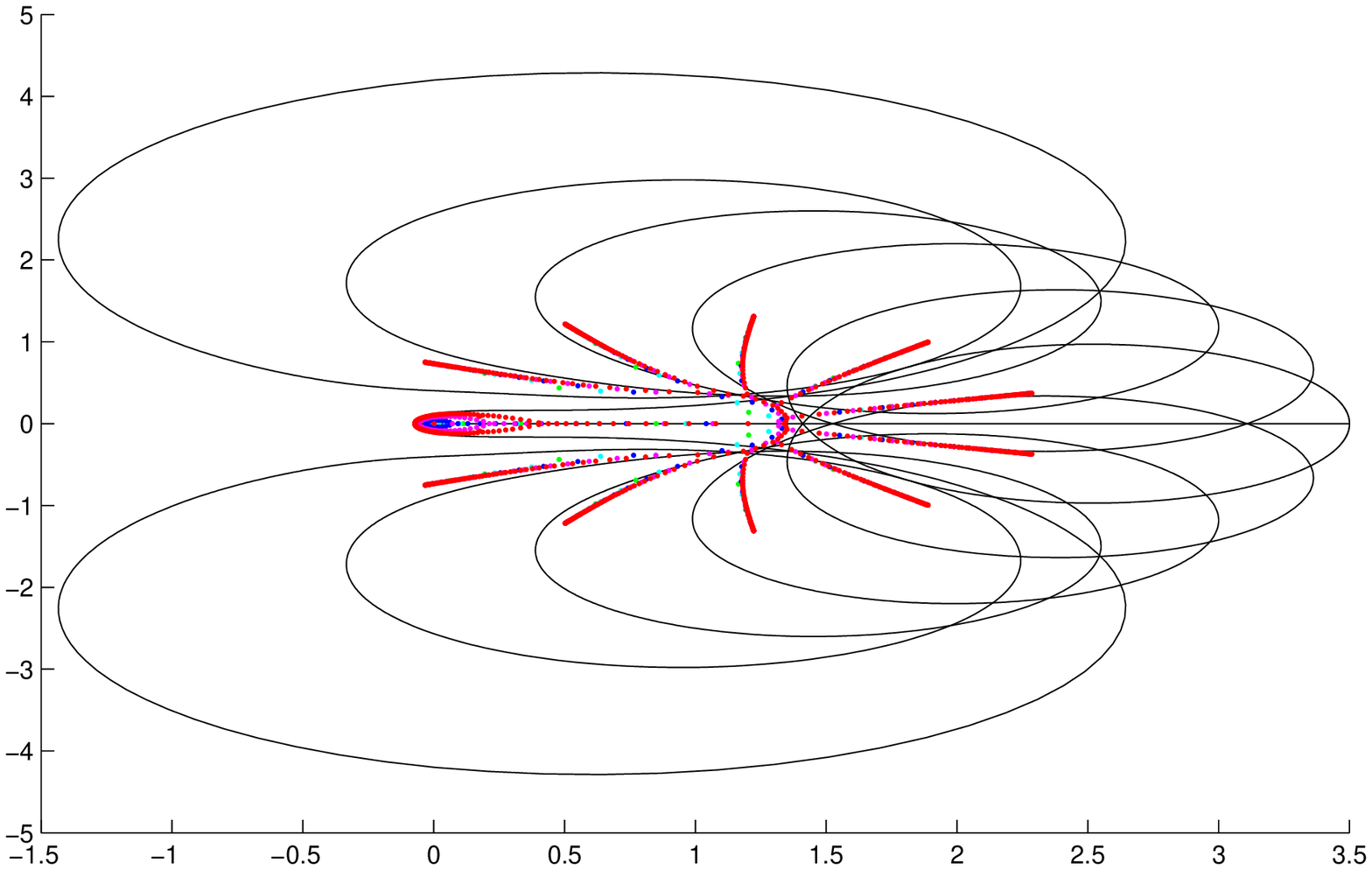} \\
\\
Figure 3. The curve $\curve(A_{\infty,k,t})$ (black) and the sets 
$\sigma(A_{50,k,t})$ (green),  \\
$\sigma(A_{100,k,t})$ (cyan),  $\sigma(A_{200,k,t})$ (blue),
 $\sigma(A_{400,k,t})$ (magenta), \\
$\sigma(A_{800,k,t})$ (red). Here $k=10$, $t=0.4$. \\ 
\end{tabular}

\begin{tabular}{c} 
\epsfysize=9.0cm       
\epsfbox{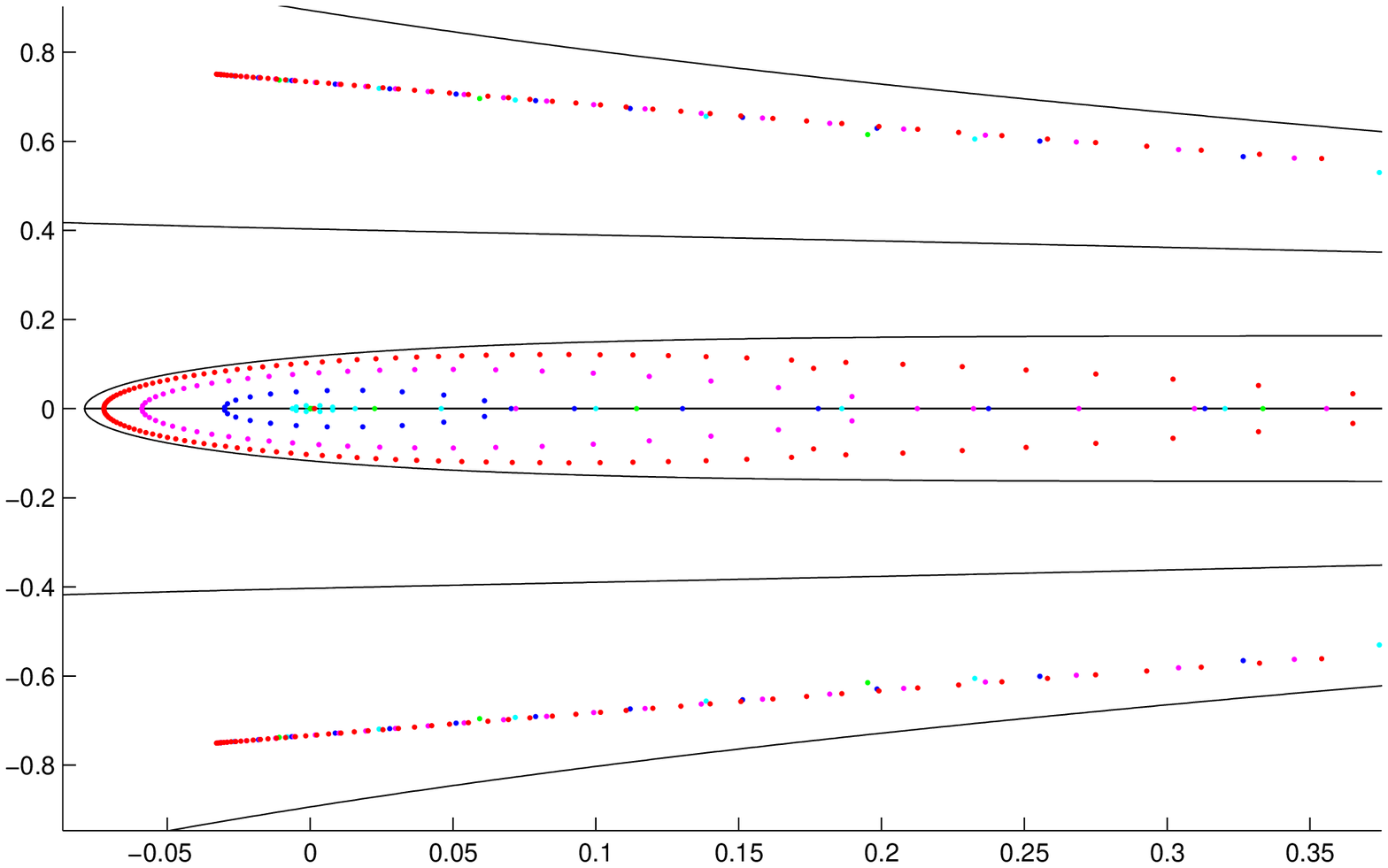} \\
\\
Figure 4. Figure 3 zoomed in the part of $\sigma(A_{\infty,k,t})$ \\
around the negative real axis.\\
\end{tabular}
\end{center}

\vskip 0.6cm

\section{Open problems} The following questions appear to deserve 
further investigation in connection with the GKK\index{GKK matrix} 
{$\tau$-matrix} \index{tau-matrix@$\tau$-matrix} problem.

\begin{description}
\item{1.} Can the matrices $A_{n,k,t}$ be approximated by $\tau$-\index{tau-matrix@$\tau$-matrix}matrices
 that are {\em strict\/} GKK\index{GKK matrix}, i.e., {$P$-matrices\index{p@$P$-matrix}} satisfying
$$ A[\alpha,\beta] A[\beta,\alpha]> 0 \qquad \forall \alpha,\beta\in
\set{n}, \quad \# \alpha=\# \beta=\# \alpha \cup \beta -1?$$
The matrices $A_{n,k,t}$ themselves are not strict GKK\index{GKK matrix!strict}\index{GKK matrix}. If the answer is no,

\item{1a.} Are strict GKK\index{GKK matrix!strict}\index{GKK matrix} matrices positive stable\index{stability}? Are strict GKK\index{GKK matrix!strict}\index{GKK matrix} {\em and\/} 
$\tau$-\index{tau-matrix@$\tau$-matrix}matrices positive stable\index{stability}?

\item{2.} Given $\alpha,\beta\subseteq\set{n}$ with $\# \alpha=
\# \beta$, call the number $\# \alpha -\# (\alpha\cap \beta)$ the 
{\em dispersal\/} of the minor $A[\alpha,\beta]$. The counterexample
from Section 1.3 shows it is not sufficient for stability\index{stability} of a $P$-matrix\index{p@$P$-matrix}
$A$ that  the inequalities 
$$ A[\alpha,\beta] A[\beta, \alpha]\geq 0 $$
hold for all minors of dispersal\index{dispersal} $\leq d\eqbd1$. 
Carlson\index{Carlson's theorem}'s theorem asserts that  the value $d=n$ is sufficient for stability\index{stability}.
What minimal value of the parameter $d$ would guarantee stability\index{stability}? 
In particular, does that value depend on $n$?
\end{description}

\noindent Also, here are two less directly related questions, which arose
in the construction of the counterexample.
\begin{description}
\item{3.} Given $n\in \N$ and positive numbers $(p_\alpha)_
{\emptyset \neq \alpha\subseteq \set{n}}$ satisfying the generalized 
Hadamard-Fisher\index{Hadamard-Fisher inequality} inequality~(\ref{HF}), when is there a matrix $A$ such that 
$A[\alpha]=p_\alpha$ for all $\alpha$?   

\item{4.}  For a matrix $A$, let 
$$ b_j\colon=\sum_{\# \alpha=j} A[{\alpha}], \qquad 
 c_j\colon={b_j \over {n \choose j}}, \qquad  j=0,\ldots,n. $$ 
When is it true that
\begin{equation} c_j^2\geq c_{j-1}c_{j+1}, \;\;\; j=1,\ldots,n-1?
 \label{newton} \end{equation}

These inequalities are known for diagonal\index{diagonal matrix} matrices 
with positive diagonal elements (e.g.,~\cite[p.51]{HLP}) 
 and go back to Newton\index{Newton's inequalities}. Since the numbers
$c_j$ are invariant under similarity, Newton's inequalities~(\ref{newton})
also hold for all diagonalizable matrices with positive real eigenvalues.
Do the GKK\index{GKK matrix} matrices satisfy~(\ref{newton})? 

\end{description}

\chapter{Inverses of special matrices}

\section{Bounded invertibility problem} The most general setup of 
the problem of this chapter is the following.

Let $\cal A$  be a collection of  matrices  such that
\begin{equation}
\inf_{A \in {\cal A}} \min \{ |z| : z\in \sigma(A) \} >0,
\qquad  \sup_{A\in {\cal A}} \|A\|<\infty \label{assump} 
\end{equation} 
for some norm $\|\cdot\|$. What conditions on $(A_j)$ 
imply 
\begin{equation} \sup_{A\in {\cal A}} \|A^{-1} \|<\infty? \label{concl} 
\end{equation}

In the easy case when the order of  matrices is bounded above, the 
conclusion~(\ref{concl})  holds without any additional hypothesis.

\nt {\bf Proposition 8.\/} {\em  Let $\cal A$ be a collection of complex
matrices satisfying~(\ref{assump}) for some norm $\| \cdot \|$ and 
such that 
\begin{equation}
\sup_{A\in {\cal A}} \order(A)<\infty. \label{boundsize}
\end{equation} 
Then (\ref{concl}) holds. }

{\em Proof.\/} \hspace{0.3cm}
Without loss of generality, one can assume that
all matrices $A\in {\cal A}$ belong to $\C^{n\times n}$ where 
$n\eqbd\sup_{A\in {\cal A}} \order(A)$.
Indeed, just replace each $A$ by
$$ \widetilde{A}\eqbd \diag(A,I_{(n-\order(A))_+}).$$
Then~(\ref{assump}) holds for the new collection
 $\widetilde{\cal A}$. Also, $\widetilde{\cal A}$ 
satisfies~(\ref{concl})
iff $\cal A$ satisfies~(\ref{concl}). Next, since all norms on 
$\C^{n\times n}$ are equivalent, one can assume 
that $\| \cdot \|$ is an operator norm\index{operator norm} subordinate to 
a norm on $\C^n$  (also denoted by $\| \cdot\|$).

So, suppose ${\cal A} \subset \C^{n\times n}$ and~(\ref{concl}) 
is violated.  Then, by the Banach-Steinhaus\index{Banach-Steinhaus' theorem}
 theorem, there exists  $v\in \C^n$
and a sequence $(A_j)$ such that 
$$\| A^{-1}_j v \| \belowrightarrow{j\to \infty} \infty.$$ 
But the sequence $(A_j)$  is totally bounded in $\C^{n\times n}$, 
hence contains a Cauchy subsequence. Without loss, it 
is $(A_j)$ itself. The limit $A\eqbd\lim_{j\to \infty} A_j$
is invertible since 
$$ \min \{ |z| : z\in \sigma(A) \}=\lim_{j\to \infty} \min \{ |z| :
z\in \sigma(A_j) \},$$
hence $\lim \|A^{-1}_j v\|=\|A^{-1}v\|<\infty$. This contradiction shows 
that (\ref{concl}) holds for any operator norm\index{operator norm}, hence
for any norm on $\C^{n\times n}$. \hfill $\Box$

Next, let us consider the case when the order of matrices $A_j$ is not
bounded above and the norm in question is the $\infty$-norm\index{p@$p$-norm!$\infty$-norm}. This question
was posed  by K. West~\cite{W} \index{West's problem}
 for positive definite\index{Hermitian matrix!positive definite} Hermitian\index{Hermitian matrix} matrices 
in connection with a problem from econometrics. As is shown in Section 2.3,
 the answer in that case is no. However, under certain additional hypotheses
the answer is yes, as is shown next.  

\section{Boundedly invertible collections of Hermitian\index{Hermitian matrix} matrices}

One of the possible restrictions on the collection $(A_j)$ that ensures 
that~(\ref{concl}) holds is (uniform) bandedness \index{banded matrix}
of matrices $A_j$,
as follows directly from a theorem of S. Demko\index{Demko's theorem}.
To recall,  a matrix $A=(a(i,j))$ is called {\em banded with band
width\/} $w$ if 
$$ a(i,j)=0 \qquad {\rm{whenever}} \quad |i-j|\geq w. $$ 

\nt {\bf Theorem (Demko\index{Demko's theorem} \cite{Dem}).\/} {\em
Let $A\in \C^{n\times n}$ be banded with band width $w$ and satisfy conditions
$\|A\|_p\leq 1$, $\|A^{-1} \|_p\leq \mu^{-1}$ for some $1\leq p \leq \infty$
and $\mu>0$. Then, with $A^{-1}\bdeq(\alpha(i,j))_{i,j=1}^n$, there are
numbers $K$ and $r\in (0,1)$ depending only on $\mu$ and $\omega$ such that
$$ |\alpha(i,j)|\leq K r^{|i-j|} \qquad \forall i,j=1,\ldots,n.$$
In particular, for any $1\leq q\leq \infty$, $\| A^{-1}\|_q\leq C$
where the bound $C$ depends only on $\mu$, and $w$.  }

Since the smallest eigenvalue of a Hermitian\index{Hermitian matrix} matrix is the reciprocal of the $2$-norm \index{p@$p$-norm!$2$-norm} of its inverse, the hypothesis of 
Demko\index{Demko's theorem}'s theorem is satisfied for $p\eqbd2$ whenever 
the collection $\cal A$ satisfies~(\ref{assump}), hence  the collection 
${\cal A}^{-1}$ is bounded in any $q$-norm ($1\leq q\leq \infty$), 
in particular, in the $\infty$-norm\index{p@$p$-norm!$\infty$-norm}.

A different restriction that ensures boundedness of $(A^{-1}_j)$,
provided that $A_j$'s are Hermitian\index{Hermitian matrix} and satisfy~(\ref{assump}), is the
oscillatory\index{oscillatory matrix} property. Recall that a matrix is 
{\em totally positive\index{totally positive matrix}\/} if
all its minors are positive. A totally nonnegative\index{totally nonnegative matrix} matrix $A$ is called 
{\em oscillatory\index{oscillatory matrix}\/} if $A^l$ is totally positive\index{totally positive matrix} for some $l\in \N$.
It is well known (see, e.g.,~\cite[p.123]{GK}) that 
$\sigma(A)$ consists of $n$ distinct positive real numbers
$\lambda_1<\cdots<\lambda_n$ and that the $k$th eigenvector $v_k$
($Av_k=\lambda_k v_k$) (unique up to a scalar multiple) has no zero 
entries and precisely $n-k$ sign changes if $A$ is oscillatory\index{oscillatory matrix}. 

\nt {\bf Theorem 2.\/} {\em Let $A\in \C^{n\times n}$ be an oscillatory\index{oscillatory matrix} Hermitian\index{Hermitian matrix}
matrix with smallest eigenvalue $\lambda_{min}$. Then 
\begin{equation} 
\|A^{-1}\|_\infty\leq {\|A\|_\infty \over \lambda^2_{min}}. \label{oscibound}
\end{equation} }

The proof will make use of the following lemma due to 
C.~de~Boor. 

\nt {\bf Lemma (de Boor\index{Boor's lemma@de Boor's lemma} \cite{dB2}).\/} {\em 
Let $A\in\C^{n\times n}$ be a nonsingular totally nonnegative\index{totally nonnegative matrix} matrix.
If, for some $x$, $y\in \C^n$, 
$$ Ax=y, \quad \sgn x(i)=\sgn y(i)=(-1)^{i-1}, $$
then $$ \|A^{-1}\|_\infty \leq {\|x\|_\infty\over\min_i |y(i)| }. $$  }

{\em Proof of Theorem 2.\/}  Since $A$ is Hermitian\index{Hermitian matrix}, its eigenvector 
$v$ corresponding to the eigenvalues $\lambda_{min}$ can be uniquely 
determined from the minimization problem
$$ v^* A v \to \min$$
once one of the entries of $v$ is fixed.
Let $v(1)=1$, so that $v=[1\;\; \widetilde{v}]^T$. Since $A$ is real,
all the entries of $v$ are necessarily real. Let $A$ be partitioned
conformably to $v$:
$$ A=\left( \begin{array}{cc} a(1,1) & A(1,2{:}n) \\
A(1,2{:}n)^T  & A(2{:}n) \end{array} \right).$$  
Then $v^T A v=a(1,1) + 2A(1,2{:}n) \widetilde{v} +\widetilde{v}^T A(2{:}n) 
\widetilde{v} $ achieves its minimum at $$ \widetilde{v}\eqbd 
-A(2{:}n)^{-1}A(1,2{:}n)^T.$$
 By the eigenvalue interlacing\index{eigenvalue interlacing} property of $A$, 
$\min \sigma(A(2{:}n))\geq \lambda_{min}$. Since $\|A(1,2{:}n)^T\|\leq 
\|A\|_1$,  this yields
$$ \| \widetilde{v} \|_\infty \leq \|\widetilde{v} \|_2 \leq 
\|A(2{:}n)^{-1}\|\cdot \|A(1,2{:}n)^T\|_2 \leq {\| A\|_\infty \over \lambda_{min}}. 
$$
The same argument can be applied to the case when any other entry of $v$ is 
set to be $1$. Since the eigenvector $v$ is unique up to multiplication
by a  scalar, this means
$$ {\|v\|_\infty \over \min_i |v(i)|} \leq {\|A\|_\infty \over \lambda_{min}}
\qquad {\rm whenever} \qquad Av=\lambda_{min} v.$$
Now apply de Boor\index{Boor's lemma@de Boor's lemma}'s lemma with $x\eqbd v$, 
$y\eqbd \lambda_{min} v$ and obtain~(\ref{oscibound}).   \hfill $\Box$

\section{On shifted Hilbert matrices and their companions}

However, the Hermitian\index{Hermitian matrix} property alone is not sufficient for the implication
{(\ref{assump}) $\Longrightarrow$ (\ref{concl})}, as is demonstrated
by the two examples below. For those counterexamples, one needs several
additional notions. 
The matrix $H_n\eqbd \left( {1\over i+j-1} \right)_{i,j=1}^n$ is known
as the {\em Hilbert matrix\/}\index{Hilbert matrix}. It has been a subject 
of extensive studies
and has served as an example of many unusual phenomena in operator theory.
M.-D. Choi in~\cite{Ch} used what he called the {\em companion of 
the Hilbert matrix\/}\index{Hilbert matrix!companion of},
viz.\ the matrix $ C_n \eqbd \left(1\over \max\{i,j\} \right)_{i,j=1}^n$.
It turns out that the matrix $C_n$ belongs to the special class of 
ultrametric\index{ultrametric} matrices\footnote{This matrix is also from
the class  of single-pair or (`one-pair') 
matrices\index{single-pair matrix} due to Gantmacher and Krein~\cite[p.113]{GK}.} introduced by Nabben and Varga~\cite{NV} as a generalization of the notion 
of strictly ultrametric\index{ultrametric!strictly} matrices attributed by 
Mart\'{\i}nez, Michon, and San Mart\'{\i}n~\cite{MMS} to 
C.~Dellacherie~\cite{De}\footnote{Ultrametricity was first introduced in
 connection with
$p$-adic number theory. Recall that a distance $d$ on a space $X$ is said to 
be ultrametric \index{ultrametric} if it satisfies the inequality
$$ d(x,y) \leq \max \{ d(x,z),d(z,y) \} \qquad \forall x,y,z\in X.$$ 
By Dellacherie's definition, a symmetric matrix $A\in \C^{n\times n}$ 
is {\em ultrametric\/} if there exists an ultrametric distance $d$ on 
$\set{n}$ such that 
$$ d(i,j)=d(i,k) \qquad \mbox{\rm iff} \qquad a(i,j)=a(i,k). $$}.  
  
A matrix $A\bdeq (a(i,j))_{i,j=1}^n$ is  {\em ultrametric\index{ultrametric}\/} if
\begin{eqnarray}
&& A=A^*, \qquad A\geq 0 \nonumber \\
&& a(i,j)\geq \min \{ a(i,k),a(k,j)\} \qquad \forall \; i,j,k\in \set{n}
 \nonumber \\
&& \mbox{\rm and} \nonumber \\
&& a(i,i) \geq \max \{ a(i,k) : k\in \set{n} {\setminus} \{i\}\} \qquad
\forall i\in \set{n}. \label{iii} 
\end{eqnarray}
If the inequality~(\ref{iii}) is strict for all $i\in \set{n}$, $A$
is called {\em strictly ultrametric\index{ultrametric!strictly}\/}. 

Finally, before stating the result of 
Mart\'{\i}nez, Michon, and San Mart\'{\i}n, 
recall that a matrix $A\bdeq (a(i,j))\in \C^{n\times n}$ is {\em row 
(column) diagonally dominant\index{diagonally dominant}\/} if
\begin{eqnarray} 
&a(i,i)\geq \sum_{j\in \set{n}\setminus\{i\}} |a(i,j)| \qquad \forall 
i\in \set{n}& \label{row} \\
&( a(i,i)\geq \sum_{j\in \set{n}\setminus\{i\}} |a(j,i)| \qquad \forall 
i\in \set{n}). & \label{column}
\end{eqnarray}  
If all the inequalities~(\ref{row}) (the inequalities~(\ref{column})) are
strict, $A$ is called {\em strictly\/} row (column) diagonally dominant\index{diagonally dominant!strictly}.\index{Mart\'{\i}nez-Michon-San Mart\'{\i}n's theorem}

\nt {\bf Theorem (Mart\'{\i}nez, Michon, and San Mart\'{\i}n~\cite{MMS}).\/} {\em
The inverse of a strictly ultrametric\index{ultrametric!strictly} matrix
is a symmetric strictly diagonally dominant $M$-matrix\index{m@$M$-matrix}. }

Now one can construct the following counterexample to the implication
{(\ref{assump}) $\Longrightarrow$ (\ref{concl})}. 

\nt {\bf Proposition 9.\/} {\em  Let $\alpha>0$ and let  
$ A_n \colon = (\alpha I_n + C_n)^{-1}$. Then the collection
$(A_n)_{n\in \N}$ of Hermitian\index{Hermitian matrix} positive definite\index{Hermitian matrix!positive definite} matrices 
satisfies~(\ref{assump}) but  does not satisfy~(\ref{concl}). }

{\em Proof.\/} Subtracting of the $j$th column from the $j-1$st column
of the matrix $C_n$, for $j=2,\ldots, n$, one verifies that $\det C_n>0$
for all $n\in \N$. Hence $C_n$ is a Hermitian\index{Hermitian matrix} positive definite\index{Hermitian matrix!positive definite} matrix.
So, $\min \sigma(\alpha I_n + C_n) >\alpha$. By~\cite[Problem V]{Ch}, 
$\|C_n\|_2 \leq 4$, so $\max \sigma (\alpha I_n +C_n)\leq 4+\alpha$.
Thus, $\sigma(A_n)\subset [1/(\alpha+4),1/\alpha]$ for any $n\in \N$.

The matrices $C_n$ are ultrametric\index{ultrametric}, hence the matrices $\alpha I_n +C_n$
are strictly ultrametric\index{ultrametric!strictly}, so by the theorem of Mart\'{\i}nez\index{Mart\'{\i}nez-Michon-San Mart\'{\i}n's theorem}, Michon and
San Mart\'{\i}n, their inverses $A_n$ are diagonally dominant\index{diagonally dominant} $M$-matrices\index{m@$M$-matrix}. 
But any diagonal entry of $A_n$ is bounded above by $\|A_n\|_2\leq 
1/\alpha$, so the $\infty$-norm\index{p@$p$-norm!$\infty$-norm} of $A_n$ is 
bounded above by $2/\alpha$. On the other hand,
$$\| A_n^{-1} \|_\infty =\|\alpha I_n + C_n \|_\infty \approx 
\alpha + \ln n \belowrightarrow{n\to \infty} 
\infty. $$ 
So, the collection $(A_n)$ satisfies~(\ref{assump}) 
but violates~(\ref{concl}). \hfill $\Box$

\nt {\bf Proposition 10.\/} {\em For large enough $\alpha >0$, the collection
$(A_n\eqbd\alpha I_n +H_n)^{-1}$ satisfies~(\ref{assump}) but does not 
satisfy~(\ref{concl}). }

{\em Proof.\/} Note that $ C_n-H_n \geq 0$ and estimate 
$\| C_n-H_n \|_\infty$. \begin{eqnarray*} 
\| C_n-H_n\|_\infty & = & \max_{i\in \set{n}}
\left( \sum_{j=1}^i \left( {1\over i} -{1 \over i+j-1}\right) + 
\sum_{j=i+1}^\infty  \left({1 \over j}-{1 \over (i+j-1)}\right) \right) \\
&= & \max_{i\in \set{n}}
\left( \sum_{j=1}^i \left( {1\over i} -{1 \over i+j-1}\right) + 
\sum_{j=i}^{2i-1} {1\over j} \right) \leq \max_{i\in\set{n}} {2i\over i}=2.
\end{eqnarray*} 

Recall that 
$$ \| A^{-1}\| \leq {\|B^{-1}\| \over 1-\|A-B\|\cdot \|B^{-1}\|} \qquad
{\rm whenever} \; \|A-B\| \cdot \| B^{-1} \|<1$$
for any operators $A$, $B$ and operator norm\index{operator norm}
 $\| \cdot \|$.
The proof of Proposition~9 demonstrated that $\|(\alpha I_n +C_n)^{-1}
\|_\infty\leq 2/\alpha$. Hence,  if $\alpha>4$, then 
$$ \|C_n-H_n\|_\infty \| (\alpha I_n+C_n)^{-1} \|_\infty 
\leq {4\over \alpha} <1,$$
hence $\|(\alpha I_n + H_n)^{-1}\|_\infty\leq {2\over \alpha -4}.$
Thus, the collection $(A_n=(\alpha I_n + H_n)^{-1})$ satisfies~(\ref{assump})
but violates~(\ref{concl}). \hfill $\Box$
  
\section{Inverses of nonnegative Hermitian\index{Hermitian matrix} Toeplitz\index{Toeplitz matrix} matrices}
 
The last example in the same spirit deals with Hermitian\index{Hermitian matrix} Toeplitz\index{Toeplitz matrix}
matrices. 

\nt {\bf Proposition 11.\/} {\em Let $T_k$ be the infinite (upper triangular) 
Toeplitz\index{Toeplitz matrix} matrix with symbol\index{Toeplitz matrix!symbol of} $\sm(T_k)(s)\colon =1+s+cs^k$, where $c$ is any 
complex number with  $|c|>2$ (chosen to be positive if the matrices $A_n$ 
must be  nonnegative).  Set  $A_k\colon=T_k T_k^*$  and let 
$A_{n,k}$ be the leading principal submatrix of $A_k$ of order $n$.
Then the collection $(A_{n,k})$ of positive definite\index{Hermitian matrix!positive definite} Hermitian\index{Hermitian matrix} Toeplitz\index{Toeplitz matrix}
matrices satisfies~(\ref{assump}) but violates~(\ref{concl}).  } 

{\em Proof.\/} \hspace{0.3cm} By the spectral theory of Toeplitz\index{Toeplitz matrix} matrices,
which was already discussed in Section~1.5, $\sigma(A_k)$ equals 
the  set of values of its symbol\index{Toeplitz matrix!symbol of} on the unit circle $|s|=1$. Notice that
 $|\sm(T_k)(s)|\geq |cs^k|-|1+s|\geq |c|-2>0$ whenever 
$|s|=1$, so $$ \inf_{k\in \N} \min \sigma(A_k)>0.$$
On the other hand, the matrices $A_k$ have at most $9$ nonzero diagonals
with the absolute value of each nonzero term at most $2$, so   
$$ \sup_{k\in \N} \|A_k\|_\infty \leq 18.$$
Finally, $A_k^{-1}=T_k^{*-1}T_k^{-1}$ and the symbol\index{Toeplitz matrix!symbol of} $\sm(T_k^{-1})$ of 
$T^{-1}_k$ is
$$ {1\over \sm(T_k)(s)} =1-s+s^2-\cdots +(-1)^{k-1} s^{k-1}+\cdots. $$
So, the $(k-1)\times (k-1)$ leading principal submatrix of
$A_k^{-1}$ has the form
$$ \left( \begin{array}{rrrrr} 
1 &  -1 & 1 &  \cdots & (-1)^{k} \\
-1 & 2 &  -2 & \cdots & (-1)^{k-1} 2  \\
1 & -2 &  3 & \cdots &  (-1)^{k-2} 3 \\
\vdots & \vdots & \vdots & \ddots & \vdots \\
(-1)^{k}  &  (-1)^{k-1} 2 & (-1)^{k-2} 3 &  \cdots &  k-1 \end{array}
\right),$$  
hence $ \sup_{k\in \N} \| A_k^{-1}  \|_1=\infty.$

Since the matrices $A_k$ are Hermitian\index{Hermitian matrix}, the limit set of the eigenvalues of 
$A_{n,k}$ (as $n$ tends to infinity) coincides with $\sigma(A_k)$ 
(see, e.g.,~\cite{BGS}). Since the (elementwise) limit of $A^{-1}_{n,k}$
is  $A_k^{-1}$ by~\cite[p.74]{GF}, the collection $(A_{n,k})$ 
satisfies~(\ref{assump}) but not~(\ref{concl}).  \hfill $\Box$

\section{Least-squares spline projection matrices}   

An interesting problem of the same type arises in spline theory.

Bounding the $\infty$-norm\index{p@$p$-norm!$\infty$-norm} of the ($L_2$-)orthogonal projector
onto splines leads to matrices of the specific form
\begin{equation} A_{n,k,{\bf t}}(i,j)\colon={k \over t_{i+k}-t_i} \int B_{ik} B_{jk}.
\label{spmatr} \end{equation}
Here $n$, $k\in \N$, ${\bf t}$ is a nondegenerate knot sequence
with $n+k$ knots, $B_{ik}$ is the $i$-th $B$-spline of order $k$ 
for the sequence ${\bf t}$, $S_{k,{\bf t}}\colon= \spa \{ B_{1k},
\ldots, B_{nk} \}$, and $Lf$ is the least squares approximation
to $f\in L_\infty[t_1,t_{n+k}]$ by elements of $S_{k,t}$. 
C.~de~Boor~\cite{dB1} showed that there exists a positive constant 
$C_k$ such that
$$ C_k \| A_{n,k,{\bf t}}^{-1}\|_\infty \leq \|L\|_\infty
\leq \| A_{n,k,{\bf t}}^{-1}\|_\infty $$
and conjectured \index{Boor's conjecture@de Boor's conjecture} that, for $k$ fixed,
\begin{equation}
\sup_{n,{\bf t}} \| A_{n,k,{\bf t}}^{-1} \|_\infty <\infty. 
\label{conj}
\end{equation}    
This conjecture was recently proved by A. Shadrin~\cite{Sha} using 
sophisticated
tools from spline theory to construct vectors $x$ and $y$ appearing in 
de Boor\index{Boor's lemma@de Boor's lemma}'s lemma with the $\min |y(i)|$ and $\max |x(i)|$ depending only on 
$k$ but not on the knot sequence ${\bf t}$ or~$n$ and conclude, using the
lemma, that~(\ref{conj}) holds.

The matrices $A_{n,k,{\bf t}}$ are known to be 
oscillatory\index{oscillatory matrix} and diagonally similar to 
Hermitian\index{Hermitian matrix} matrices, i.e.,
$$ A_{n,k,{\bf t}}=D_{n,k,{\bf t}} H_{n,k,{\bf t}} D^{-1}_{n,k,{\bf t}}$$
where $D_{n,k,{\bf t}}$ are diagonal\index{diagonal matrix} with positive diagonal entries. (For sure, the condition number\index{condition number}
 of the matrices $D_{n,k,{\bf t}}$ is  not bounded above.)  Moreover, 
$A_{n,k,{\bf t}}$ are banded  \index{banded matrix}
with band width $k$. Finally, the smallest
eigenvalue of $H_{n,k,{\bf t}}$ is known to be bounded away from zero
independently of ${\bf t}$ and $n$, so the collection $(A_{n,k,{\bf t}})$
satisfies~(\ref{assump}). (All those facts can be found in, 
e.g.,~\cite[p.401--406]{DL}.)
However, the above conditions are not sufficient for~(\ref{concl}) (that 
is,~(\ref{conj})) to hold. In particular, an application of Theorem~2
to the eigenvector $v_{min}$ of $H_{n,k,{\bf t}}$ corresponding to its 
smallest eigenvalue demonstrates that the eigenvector $D_{n,k,{\bf t}}v$
of $A_{n,k,{\bf t}}$ corresponding to the same eigenvalue cannot be 
used to prove~(\ref{conj}). 

Hence the following (somewhat vaguely formulated) problem:

\nt {\bf Problem.\/} {\em
Single out an additional property of the matrices $A_{n,k,{\bf t}}$ 
to obtain a simple matrix theoretic proof of~(\ref{conj}).  }

\printindex


\begin{thebibliography}{111}

\bibitem{An} T. Ando\index{Ando's theorem},  Totally positive\index{totally positive matrix} matrices, 
Linear Algebra Appl. 90 (1987) 165--219.

\bibitem{BP} A. Berman, R. Plemmons, Nonnegative matrices in the mathematical
sciences. SIAM, Philadelphia, 1994.

\bibitem{Bh} R. Bhatia, Matrix analysis. Springer-Verlag, New York, 1997.

\bibitem{Bi} M. Biernacki\index{Biernacki's theorem}, Sur les \'{e}quations alg\'{e}briques contenant
des param\'{e}tres arbitraires, Bull. Acad. Polon. Sci. (1927) 541--685.

\bibitem{dB1} C. de Boor, \index{Boor's conjecture@de Boor's conjecture} The quasi-interpolant as a tool in elementary
polynomial spline theory, {\em in\/} ``Approximation Theory''
 (C. G. Lorentz, ed.), pp. 269--276, Academic Press, New York, 1973.

\bibitem{dB2} C. de Boor\index{Boor's lemma@de Boor's lemma},  On a max-norm bound for the least-squares spline
approximant, {\em in\/} ``Approximation and Function Spaces,'' Proceedings,
International Conference, Gdansk, 1979, pp.~163--175, New York, 1981.

\bibitem{BGS} A. B\"{o}ttcher, S. M. Grudsky, B. Silbermann, Norms of 
inverses, spectra, and pseudospectra of large truncated Wiener-Hopf operators
and Toeplitz\index{Toeplitz matrix} matrices, New York J. Math. 3 (1997), 
1--31 (electronic).    

\bibitem{C1} D. Carlson\index{Gantmacher-Krein-Carlson's theorem},  weakly sign-symmetric\index{sign-symmetric!weakly} matrices and some
determinantal inequalities, Colloq. Math. 17 (1967) 123--129.

\bibitem{C2} D. Carlson\index{Carlson's theorem},  A class of positive stable\index{stability} matrices,
J. Res. Nat. Bur. Standards Sect. B 78 (1974) 1--2.

\bibitem{Ca} A. Cauchy, Sur l'equation a l'aide de laquelle on d\'{e}termine
les inegalit\'{e}s seculaires de mouvement de plan\`{e}tes, Oeuv. comp.
9(2) (1829), 174--195. 

\bibitem{Ch} M.-D. Choi, Tricks or treats with the Hilbert matrix, Amer. Math.
Monthly 90 (1983), no.5, 301--312.

\bibitem{D} K. M. Day\index{Day's theorem}, Toeplitz\index{Toeplitz matrix} matrices generated by the Laurent\index{Laurent series} series 
expansion of an arbitrary rational function, Trans. Amer. Math. Soc.
206 (1973) 224--245.

\bibitem{De} C. Dellacherie, Private communication (1985).

\bibitem{Dem} S. Demko\index{Demko's theorem}, Inverses of band matrices and local convergence of 
spline projections, SIAM J. Numer. Anal. 14 (1977), no.4, 616--619.

\bibitem{DL} R. DeVore and G. Lorentz, Constructive approximation.
Springer-Verlag, Berlin, 1993.

\bibitem{ES} G. M. Engel\index{Engel-Schneider's conjecture} and H. Schneider,  The Hadamard-Fisher\index{Hadamard-Fisher inequality} 
inequality for a class of matrices defined by eigenvalue\index{eigenvalue monotonicity} monotonicity,
Linear and Multilinear Algebra  4 (1976) 155--176.

\bibitem{GK} F.R. Gantmacher, M.G. Krein, Oscillation matrices and kernels
and small vibrations of mechanical systems. Gostechizdat, 1950.
\index{Gantmacher-Krein-Carlson's theorem}

\bibitem{GF} I. C. Gohberg and I. A. Fel'dman, Convolution equations and
projection methods for their solution. Translations of Mathematical Monographs,
vol. 41, AMS, Providence,  1974. 

\bibitem{Fr} S. Friedland, Weak interlacing\index{eigenvalue interlacing} property of totally positive\index{totally positive matrix} 
matrices, Linear Algebra Appl. 71 (1985) 95--100. 

\bibitem{GS} I. Gohberg and A. A. Semencul, On inversion of finite-section
Toeplitz\index{Toeplitz matrix} matrices and their continuous analogues (in Russian),  Matem. Issled., Kishinev, 7(2) (1972) 201--224. 

\bibitem{HLP} G. Hardy, J. E. Littlewood, and G. Polya, Inequalities, 
Cambridge Univ. Press, Cambridge-New York, 1952.

\bibitem{HR} G. Heinig and K. Rost, Algebraic methods for Toeplitz\index{Toeplitz matrix}-like matrices and operators, Birkh\"{a}user Verlag, Basel-Boston,
Mass., 1984.

\bibitem{He} D. Hershkowitz\index{Hershkowitz' conjecture},  Recent directions in matrix stability\index{stability},
Linear Algebra Appl. 171 (1992) 161--186.

\bibitem{HB} D. Hershkowitz and A. Berman,  Notes on $\omega$- and 
$\tau$-matrices, Linear Algebra Appl. 58 (1984) 169--183.

\bibitem{Ho} O. Holtz, Not all GKK\index{GKK matrix} 
$\tau$-\index{tau-matrix@$\tau$-matrix}matrices are stable.
Linear Algebra Appl. 291 (1999) 235--244.

\bibitem{Ke} R. B. Kellogg\index{Kellogg's theorem}, On complex eigenvalues 
of $M$ and $P$  matrices, Numer. Math. 19 (1972) 170--175. 

\bibitem{MMS} S. Mart\'{\i}nez\index{Mart\'{\i}nez-Michon-San Mart\'{\i}n's theorem}, G. Michon, and J. San Mart\'{\i}n, Inverses of
ultrametric \index{ultrametric} matrices are of Stieltjes type,
SIAM J. Matrix Anal. Appl. 15 (1994) 98--106.

\bibitem{M} V. Mehrmann,  On some conjectures on the spectra of
$\tau$-\index{tau-matrix@$\tau$-matrix}matrices, Linear and Multilinear Algebra  16 (1984) 101--112.

\bibitem{NV} R. Nabben and R. Varga, Generalized ultrametric matrices
\index{ultrametric} -- a class of inverse 
$M$-matrices\index{m@$M$-matrix}, Linear Alg. Appl. 220 (1995) 365--390.

\bibitem{ScS} P. Schmidt and F. Spitzer, The Toeplitz\index{Toeplitz matrix} matrices of an arbitrary
Laurent\index{Laurent series} polynomial, Math. Scand. 8 (1960), 15--38.

\bibitem{Sha} A. Shadrin, The $L_\infty$-norm of the $L_2$-spline-projector
is bounded independently of a knot sequence. A proof of 
de~Boor\index{Boor's conjecture@de Boor's conjecture}'s conjecture,
manuscript (1999).

\bibitem{TT} O. Taussky, Research problem, Bull. Amer. Math. Soc. 64 (1958) 
124.

\bibitem{V} R. Varga\index{Varga's conjecture},  Recent results in linear algebra and its
applications (in Russian), in: Numerical Methods in Linear Algebra,
Proceedings of the Third Seminar of Numerical Applied Mathematics, 
Akad. Nauk SSSR Sibirsk. Otdel. Vychisl. Tsentr, 
Novosibirsk, 1978, pp. 5-15.

\bibitem{W} K. West, Private communication (1999).

\end{thebibliography}
\end{document}